\documentclass[ijoo,nonblindrev]{informs-ijoo}
\usepackage{natbib}
 \bibpunct[, ]{(}{)}{,}{a}{}{,}%
 \def\bibfont{\small}%

\TheoremsNumberedThrough
\ECRepeatTheorems
\EquationsNumberedThrough \MANUSCRIPTNO{}

\usepackage{appendix}

\usepackage[font=small]{subcaption}
 \usepackage{graphicx}

 \usepackage{nomencl} %
\usepackage{pgfplots}
\usepgfplotslibrary{fillbetween}
\pgfplotsset{
  compat=newest,
  xlabel near ticks,
  ylabel near ticks
}

\usepackage{soul} 

 \usetikzlibrary{calc} 
 \usepackage{circuitikz} 
 \usepackage{bm}
\usepackage{nomencl} 
 \makenomenclature
 
\usepackage{todonotes}

\newcommand{\bx}{\mathbf{x}}
\newcommand{\by}{\mathbf{y}}

\newcommand{\bz}{\mathbf{z}}
\newcommand{\bzp}{\mathbf{z}^+}
\newcommand{\bzn}{\mathbf{z}^-}
\newcommand{\bzpA}{\mathbf{z}^{A-}}
\newcommand{\bznA}{\mathbf{z}^{A-}}
\newcommand{\br}{\mathbf{r}}
\newcommand{\bbee}{\mathbf{e}}
\newcommand{\bv}{\mathbf{v}}
\newcommand{\bh}{\mathbf{h}}

\newcommand{\balpha}{\bm{\alpha}}
\newcommand{\byhat}{\hat{\mathbf{y}}}
\newcommand{\bshat}{\hat{\mathbf{s}}}
\newcommand{\bA}{\mathbf{A}}
\newcommand{\bB}{\mathbf{B}}
\newcommand{\bg}{\mathbf{g}}
\newcommand{\bzero}{\mathbf{0}}
\newcommand{\bone}{\mathbf{1}}
\newcommand{\cU}{\mathcal{U}}
\newcommand{\cY}{\mathcal{Y}}
\newcommand{\cH}{\mathcal{H}}

\newcommand{\bty}{\mathbf{\tilde{y}}}
\newcommand{\bts}{\mathbf{\tilde{s}}}
\newcommand{\buy}{\mathbf{\underline{y}}}
\newcommand{\boy}{\mathbf{\overline{y}}}
\newcommand{\bus}{\mathbf{\underline{s}}}
\newcommand{\bos}{\mathbf{\overline{s}}}

\newcommand{\bc}{\mathbf{c}}
\newcommand{\bSigma}{\mathbf{R}}

\newcommand{\bxl}{\bx^{\ell}}
\newcommand{\byl}{\by^{\ell}}
\newcommand{\bzpl}{\bz^{+\ell}}
\newcommand{\bznl}{\bz^{-\ell}}
\newcommand{\bxz}{\bx^{0}}
\newcommand{\byz}{\by^{0}}
\newcommand{\bzzp}{\bz^{+0}}
\newcommand{\bzzn}{\bz^{-0}}

\newcommand{\Gammal}{\Gamma^{\ell}}
\newcommand{\Gammaz}{\Gamma^{0}}

\newcommand{\M}{[\mathbf{N}]}

\newcommand{\MB}{[\mathbf{R}(\Gamma)]}
\newcommand{\MBm}{[\mathbf{R}(m)]}

\newcommand{\MAm}{[\mathbf{R}^{A}(m)]}
\newcommand{\MBl}{[\mathbf{R}(\Gammal)]}
\newcommand{\MBz}{[\mathbf{R}(\Gammaz)]}
\newcommand{\MAA}{[\mathbf{R}^{A}(\Gamma^A)]}

\newcommand{\suj}{\underline{s}_j}
\newcommand{\soj}{\overline{s}_j}
\newcommand{\shj}{\hat{s}_j}

\newcommand{\yuj}{\underline{y}_j}
\newcommand{\yoj}{\overline{y}_j}
\newcommand{\yhj}{\hat{y}_j}

\newcommand{\sign}{\mathrm{sgn}}

\definecolor{myGreen}{RGB}{0,150,50}

\usepackage{interval}

\begin{document}
\RUNAUTHOR{Dehghani Filabadi and Mahmoudzadeh}
\RUNTITLE{Effective Budget of Uncertainty for Classes of Robust Optimization} 
\TITLE{Effective Budget of Uncertainty for \\ Classes of Robust Optimization} 
\ARTICLEAUTHORS{%
\AUTHOR{Milad Dehghani Filabadi}
\AFF{Department of Management Sciences, University of Waterloo, \EMAIL{m23dehgh@uwaterloo.ca}} 
\AUTHOR{Houra Mahmoudzadeh}
\AFF{Department of Management Sciences, University of Waterloo, \EMAIL{houra.mahmoudzadeh@uwaterloo.ca}}
} 

\ABSTRACT{%
Robust optimization (RO) tackles data uncertainty by optimizing for the worst-case scenario of an uncertain parameter and, in its basic form, is sometimes criticized for producing overly-conservative solutions. To reduce the level of conservatism in RO, one can use the well-known budget-of-uncertainty approach which limits the amount of uncertainty to be considered in the model. In this paper, we study a class of problems with resource uncertainty and propose a robust optimization methodology that produces solutions that are even less conservative than the conventional budget-of-uncertainty approach. We propose a new tractable two-stage robust optimization approach that identifies the ``ineffective'' parts of the uncertainty set and optimizes for the ``effective'' worst-case scenario only. In the first stage, we identify the effective range of the uncertain parameter, and in the second stage, we provide a formulation that eliminates the unnecessary protection for the ineffective parts, and hence, produces less conservative solutions and provides intuitive insights on the trade-off between robustness and solution conservatism. We demonstrate the applicability of the proposed approach using a power dispatch optimization problem with wind uncertainty. We also provide examples of other application areas that would benefit from the proposed approach.
}%
\KEYWORDS{robust optimization; budget of uncertainty; power dispatch optimization}
\maketitle
\section{Introduction} 
Contrary to deterministic optimization which assumes the parameters of an optimization problem are known with certainty, robust optimization (RO) \citep{soyster1973convex,el1997robust,ben1998robust,ben2002robust,bertsimas2003robust} takes into account that such perfect information is not always available due to reasons such as measurement errors, round-off computational errors, and even forecasting inaccuracies. 
An RO model often assumes that an uncertain parameter can take any value within a given uncertainty set and finds the optimal solution under the worst-case possible scenario within that set. Examples of applications of RO are widespread in areas such as portfolio selection \citep{hassanzadeh2014robust}, network flows \citep{atamturk2007two}, inventory management \citep{ang2012robust}, cancer treatment \citep{chan2014robust}, and power dispatch problems \citep{li2015robust,li2016adaptive}. We refer the reader to \cite{bertsimas2011theory} and \cite{gabrel2014recent} for comprehensive reviews on the applications of robust optimization.  

{ In addition to RO, there are other methods for handling uncertainty in optimization problems in the literature. For instance, stochastic optimization (SO) \citep{ermoliev1988numerical,schneider2007stochastic} considers a probability distribution on the uncertain parameter(s), and 
uncertainty quantification (UQ) \citep{smith2013uncertainty} uses statistical tools such as Monte-Carlo sampling or stochastic collocation to model uncertainties. However, these methods are often computationally challenging \citep{chatterjee2017improved}, their results are sensitive to the statistical information on the data \citep{ben2009robust}, and may not be sufficiently robust in certain applications \citep{imholz2020robust}. 
Such statistical information on the data has also been combined with RO in a distributionally robust optimization (DRO) methodology \citep{delage2010distributionally} which reduces the sensitivity of the solutions to the statistical data, but still requires prior information on the uncertain parameters, which may not always be available. 

In recent literature, traditional RO models are sometimes criticized for producing overly-conservative solutions since they optimize for the worst-case scenario of the uncertain parameter(s). Several studies have focused on reducing the conservatism of RO solutions. For example, adjustable/adaptive RO (ARO) considers a multi-stage problem where the uncertainty set and the corresponding robust solution is updated every time an observation of uncertainty is made \citep{ben2004adjustable,delage2015robust}.
Such ARO models are often computationally intractable and studies have aimed to improve the tractability of these models under certain conditions \citep{bertsimas2010finite,bertsimas2012power,iancu2013supermodularity,ardestani2016robust}. For a class of large-scale distributionally robust optimization problems, \cite{gupta2019near}  
proposed an approach to identify and remove unnecessary protection and reduce the level of conservatism. In another stream of research, \citet{ide2013relation,crespi2017quasiconvexity} and \citet{chan2017stability} studied the stability of RO problems with respect to perturbations in the size of the uncertainty. These studies emphasize the importance of reducing conservatism in different types of RO models.

Most relevant to our work,  \cite{bertsimas2004price} proposed a budget-of-uncertainty approach to reduce the level of conservatism of RO models. This budget of uncertainty is a constant 
controls how much the uncertain parameters 
can deviate from their nominal values. Particularly, a zero budget corresponds to the deterministic problem with no uncertainty, and a larger budget corresponds to a higher level of uncertainty. Ultimately, a full budget refers to complete protection against uncertainty, 
which is equivalent to the traditional worst-case approach. 
This trade-off between the level of uncertainty and its effect on the objective function of the problem is referred to as the ``price of robustness''. The budget-of-uncertainty approach was originally proposed for  ``row-wise'' polyhedral uncertainty, i.e., where the rows of the constraint matrix belong to a given set.
On the other hand, “column-wise” uncertainty, i.e., where the columns of the constraint matrix belong to a given set, are often still seen as overly conservative \citep{ouorou2016robust, minoux2008robust, minoux2012two}.}

In this paper, we show that the conventional budget-of-uncertainty approach proposed by \cite{bertsimas2004price}  
may still be overprotective for column-wise uncertainty on the right-hand side (RHS) parameters. That is, for a class of problems with uncertainty on the amount of available resources considered on the RHS of constraints, this approach may be over-protective. We demonstrate that even though this approach adjusts the amount of uncertainty considered in the model, this adjustment may not affect the solution for levels of uncertainty higher than a certain threshold. That is, protecting against the absolute worst-case is not beneficial for the system since after a certain level of uncertainty, the solution will no longer change. Hence, the additional budget considered is not \emph{effective}.

To address this issue, we define an \emph{effective worst-case scenario} and propose a formulation that protects against the effective worst case instead of the absolute worst case. We define the effective worst case as the worst scenario within a level of uncertainty that can possibly affect the optimal solution. 
We show that our proposed model is less conservative than the original budget-of-uncertainty model and empirically show that it still preserves solution robustness. The proposed approach will allow the users to explicitly observe the trade-off between the level of uncertainty and solution robustness (conservatism) in this class of problems. 
Hence, decision makers can more intuitively optimize the system for various levels of the effective budget. 
We motivate and demonstrate our approach using a power dispatch optimization problem with uncertainty in the amount of available wind power in the system. 
The specific contributions of this paper are as follows:

\begin{itemize}
\item We propose a modified 
budget-of-uncertainty approach for a class of problems with resource 
uncertainty that have ineffective budgets. We show the proposed approach is less conservative than the conventional budget-of-uncertainty approach.
\item We propose a two-stage robust optimization method to find the effective worst case which allows the user to 
explicitly control the level of solution conservatism. 
\item { We motivate the problem by providing several application areas and demonstrate the applicability of the proposed approach in power dispatch optimization problems with  
uncertainty in the amount of wind power. }
\end{itemize}

The rest of this paper is organized as follows: { Section~\ref{motivation} motivates the problem using three example applications with RHS uncertainty. 
 In Section~\ref{robust model}, the proposed two-stage robust optimization model is presented and discussed in details. Section~\ref{results} provides numerical results and analyses for the power dispatch optimization problem. Finally, concluding remarks are provided in Section~\ref{conclusion}.
}
\section{Motivation} \label{motivation}
{ In this section, we first motivate the proposed approach using a power dispatch optimization problem which is the focus of the numerical results of this work. We then provide two additional applications for which we will be providing small numerical examples in Appendix~\ref{appendix other applications}.

\subsection{Power Dispatch under Wind Uncertainty}
}
Modern power systems rely 
on the integration of low-carbon renewable resources to meet the demand for electricity. Among these renewable resources, wind energy is of special importance due to increasing penetration into power systems in the recent years so that the U.S. energy's plan for 2030, wind power is expected to have a significant contribution (20\%) in providing electricity \citep{lopez2012us}.
However, due to the inherent uncertainty in wind, it is not possible 
to forecast the exact amount of available wind power in 
the existing day-ahead electricity market. If the power system is planned based on the day-ahead wind power prediction, even small prediction errors and uncertainties associated with the amount of available wind power can make the planned power dispatch infeasible by violating the operational limits of the system, and therefore, can potentially 
lead to security and reliability issues in the power system in real-time 
\citep{lorca2015adaptive}. 

Reliability and security concerns have been the focus of power systems since the 1960s \citep{billinton1968transmission}. Particularly, the security-constrained economic dispatch (SCED) problem is concerned with operational security and reliability of the system to mitigate the risk of a system failure under unforeseen contingencies \citep{frank2016introduction}. In SCED, the goal is to find the most economical power dispatch plan while considering the operational constraints of the system (e.g., power balance, generation and ramp limits, reserve requirements, and power flow transmission constraints). 
When there is an excessive amount of wind power that cannot be absorbed by the system, the wind power can be ``curtailed'' by shutting down some or all of the wind turbines in order to maintain the system's operational security.  
However, wind curtailment is not desirable from an economical point of view since wind power is a free renewable resource while the alternatives mainly rely on costly fossil fuels. Therefore, the objective function of SCED models typically contains a high penalty for wind curtailment to encourage the power system to utilize as much wind power as possible without violating the security requirements of the system \citep{li2016adaptive}. 

Given wind power uncertainty, the SCED problem can be categorized as an optimization problem with resource 
uncertainty where the optimal wind power generation has an uncertain upper bound. Due to the high curtailment penalty, the worst-case scenario happens when the realized wind power is much higher than what the system can absorb, which in turn, results in a high total cost. 
Intuitively, a higher value selected for the budget of uncertainty corresponds to more potential available wind power. However, when the available wind is much higher than the predicted (nominal) value, 
the total available wind power may never be entirely utilized due to the operational limits of the system. Thus, 
increasing the budget of uncertainty will not change the optimal solution and would merely result in additional curtailment cost. Hence, the budget is ``ineffective'' and 
there is no change in the solution so as to mitigate the absolute worst-case scenario. This behavior is the motivation behind our proposed approach 
which optimizes the system for the ``effective'' worst-case scenario instead and produces a less conservative solution.   
We will elaborate on the SCED application in more details in the numerical results in Section~\ref{results}. 

\subsection{Other Applications} \label{sec other applications}
{ 
The proposed methodology was originally motivated by the power dispatch optimization problem, but its applicability goes far beyond it. In general, many problems with RHS budget uncertainty could benefit from the proposed approach. In this section, we provide two examples of applications for which this methodology can be applied to. 

\paragraph{Patient Scheduling: }
Consider a patient scheduling system with different priorities of patients. The goal is to optimize and allocate capacity to different priorities of patients where the number of patients in each priority is uncertain in future periods. The cost of capacity is often considered as piece-wise linear convex which means that after certain thresholds, increasing capacity becomes more and more expensive. In such a problem, the number of patient arrivals per day for each priority is considered in the RHS of service constraints and patient waiting time is associated with a penalty in the objective function \citep{mahmoudzadeh2020robust, Shahraki2020advance}. When a budget of uncertainty is included in the model, after a certain threshold, increasing the level of uncertainty will not affect the optimal solution and the additional budget would be considered as inefficient. Our proposed model can instead identify the effective worst-case scenario that can help managers in better understand and plan for it in order to minimize overall patient wait times. In Appendix~\ref{appendix patient scheduling}, we provide a small numerical example of this application for interested readers. 

\paragraph{Inventory Planning: }
In inventory problems, the customer demand is often considered as an uncertain parameter \cite{bertsimas2006robust}. In certain settings, this demand can be considered on the RHS of constraints and there is a cost for shortage or delayed demand. The capacity for holding inventory is often limited and increasing inventory more than a certain level can be extremely costly. Hence, when holding/capacity costs exceed shortage costs, any increased budget of uncertainty on demand will not have any effect on the optimal solution of the model and the identified worst-case will not be effective. Our proposed model will be able to identify the effective worst-case scenario and help in making more economical production and inventory decisions. A small numerical example of this problem is provided in Appendix~\ref{appendix inventory}. 
}

\section{The Proposed Robust Optimization Approach}\label{robust model}
In this section, we first present a nominal problem with no uncertainty that represents the type of problems we consider in this work. We then develop a basic robust optimization problem with budget of uncertainty 
to mathematically describe and motivate the proposed methodology. Lastly, we define \emph{admissible} and \emph{effective} sets and develop a new two-stage robust model for dealing with resource 
uncertainty in this class of problems. 

Throughout the rest of this paper, all inequalities are element-wise. We use bold numbers $\bzero$ and $\bone$ to denote vectors of all zeros and all ones, respectively. We use operator~$\odot$ to denote the Hadamard Product \citep{horn1990hadamard} for element-wise multiplication of vectors with equal dimensions. 
\subsection{Problem Description and Mathematical Motivation}
To set up the problem, in this section, we first present a nominal problem with no uncertainty and then build a robust version of the model with budget of uncertainty. Using the robust model, we explain the challenge of using the conventional budget-of-uncertainty approach and motivate the proposed methodology. 

\subsubsection{A Nominal Model: }
Consider formulation $\M$ as the nominal model (i.e., with no uncertainty) as follows:
\begin{subequations}\label{M2}
    \begin{align} 
\M: \quad    \min_{\bx,\by} \quad & \bc_1'\bx+\bc_2'(\byhat-\by), &\label{M2-a} \\
\quad \text{s.t.} \quad & \bA\bx + \bB\by \le \bg,  &\label{M2-b}\\
& \by \le \byhat,  &\label{M2-d}\\ 
& \bx, \by \ge \bzero. & \label{M2-e}
\end{align}
\end{subequations}
{  Problem $\M$ is a deterministic optimization problem with decision variable $\bx$ and $\by$ that are column vectors of size $p$ and $m$, respectively,} and correspond to different types of resources. The objective function consists of a linear cost of using resource 1, i.e., $\bc_1'\bx$, and a penalty for under-utilizing resource 2 compared to the available amount $\byhat$, i.e., $\bc_2'(\byhat-\by)$. The objective function maximizes the amount of resource 2 used while 
its available amount is capped in the right-hand side of constraint \eqref{M2-d}. 
Constraint~\eqref{M2-b} captures the limitations of the system 
on both $\bx$ and $\by$, where matrices $\bA$ and $\bB$ are of size $n \times p$ and $n \times m$, respectively, and $\bg$ is a vector of parameters of size $n \times 1$. 

This formulation represents a model 
which can be used for various real-world applications such as the {  multi-period} SCED problem mentioned before, 
in which the vectors $\bx$ and $\by$ correspond to the power generation of conventional generators and wind power plants, respectively. Thus, $\bc_2'(\byhat-\by)$ corresponds to the cost of wind power curtailment while considering the 
known available wind power in the RHS of constraint \eqref{M2-d}.  In what follows, we list a few assumptions that are required for the proposed methodology. These assumptions hold in the SCED model and are of practical relevance in other similar applications. For ease of mathematical notations, we will develop the proposed methodology with the general form introduced in formulation~$\M$.  We later provide complete details of the specific SCED formulation in Appendix~\ref{appendix SCED} along with an example of such a system in Appendix~\ref{rts24}.

{
Let $c_{1,(i)}$ be the $i^{\text{th}}$ element of vector $\bc_1$. Define the set of indices $\mathcal{P}=\{i : c_{2,(i)} \ge \max_k \hspace{0.1cm} c_{1,(k)} \}$.

\begin{assumption}
 \label{assum:c}
$|\mathcal{P}| \ge 1$.
\end{assumption}
Assumption~\ref{assum:c} means that at least one element in $\by$ has a corresponding under-utilization penalty cost higher than the maximum production cost of resource 1, i.e., $\bx$. In the context of SCED, for example, it means that the wind curtailment cost (i.e., the cost of shutting down the wind turbines to stop absorbing the available wind power) in at least one turbine is larger than the maximum cost of power generation using conventional generators. This assumption ensures that the optimization problem will attempt to use as much of resource 2 in at least one of the turbines, to avoid the high penalty cost. In the absence of this assumption, the wind power integration would not be beneficial.}

\begin{assumption}\label{assum:B}
In formulation~$\M$, $\bB \ge \bzero$, i.e., all elements of the matrix $\bB$ are non-negative. 
\end{assumption}
We note that Assumption~\ref{assum:B} is  a practically-relevant assumption since $\bB\by \ge 0 $ in constraint \eqref{M2-b} corresponds to the system capacity used by resource 2 (for inequality constraints) or the contribution of resource 2 in satisfying a demand (for equality constraints). Furthermore, the case of $\bB=\bzero$ corresponds to the upper/lower limit of vector $\bx$, depending on the sign of the coefficients in the $\bA$ matrix. 

\begin{assumption}\label{assum:cc}
 $\bc_1, \bc_2 \geq \bzero$, i.e., the cost vectors are element-wise non-negative.
\end{assumption}

Assumption~\ref{assum:cc} ensures that the objective function minimizes $\bx$ and maximizes $\by$. In the SCED problem, these costs vectors represent monetary values which are inherently non-negative.

\subsubsection{A Robust Model with 
Budget of Uncertainty: }

Next, we consider uncertainty in the available amount of resource 2  (e.g., the available wind in the SCED problem). 
Applying the traditional robust optimization approach of \citet{soyster1973convex} 
to formulation $\M$ would lead to an over-conservative solution since it would simply find the worst-case of each uncertain 
parameter according to the direction of the objective function independently, 
and hence, would recover a strictly conservative solution. Therefore, we  
incorporate a budget-of-uncertainty approach instead, so as to obtain a model that can control the degree of conservatism~\citep{bertsimas2004price}. 
We later make modifications to the budget-of-uncertainty approach in order to further reduce the level of conservatism by protecting against the effective worst-case scenario only. 

Let $\bty$ and $\byhat$ be the uncertain parameter and nominal value for the available amount of resource 2, respectively.  
{ For simplicity, { and in accordance with the conventional budget-of-uncertainty approach,} we start by considering a box uncertainty for $\bty$ and will generalize to other types of uncertainty in Section~\ref{sec:extension}.  
Assume that each element of the vector $\bty$ can deviate from the nominal value $\byhat$ within the interval $[\buy,\,\boy]$. We denote the underlying uncertainty set (without budget) as $\cY$ defined as  
\[\cY = \{\,\bty \mid \, \buy \leq \bty \leq \boy\, \}.\] }
Considering the budget-of-uncertainty approach, the uncertain value of $\bty$ can be written as 
\[ \bty = \byhat + \bzn \odot (\buy - \byhat) + \bzp \odot (\boy - \byhat), \]
where the vectors $\bzero \le \bzn, \bzp \le \bone$ capture the negative and positive scaled deviations from the nominal value, respectively. The total amount of these deviations is controlled by a budget parameter $\Gamma$, which is set by the planner to control the level of conservatism. This parameter dictates how much deviation from the nominal values should be considered in total.

The uncertainty set on the vector $\bty$ can be considered  
as follows:
\begin{equation}
\label{U set with conventional budget}
\cU = \Big\{ \hspace{0.051cm} \bty \in \mathbb{R}^{{m}} : \hspace{0.2cm}  \bty = \byhat + \bzn \odot (\buy - \byhat) + \bzp \odot (\boy - \byhat),
\hspace{0.2cm} \sum_{j=1}^{m} (z^-_{j}+z^+_{j}) \le \Gamma, \hspace{0.2cm} \bzero \le \bzn,\bzp \le \bone \hspace{0.051cm} \Big\},
\end{equation}

Incorporating this uncertainty set $\cU$ into the nominal model~$\M$ we derive a robust formulation with budget $\Gamma$ that minimizes the worst-case objective function value under uncertainty as follows.
\begin{subequations}\label{M2 with Budget}
\begin{align}
\MB: \quad \min_{\bx,\by} \quad & \Big\{ \mathbf{c}_1'\bx+ \max_{\bz} \mathbf{c}_2'\Big(
\byhat + \bzn \odot (\buy - \byhat) + \bzp \odot (\boy - \byhat)-\by\Big) \Big\}, & \label{M2-a budget} \\
\text{s.t} \quad & \bA\bx+\bB\by \le \bg, & \label{M2-b budget}\\ & \by \le \byhat + \bzn \odot (\buy - \byhat) + \bzp \odot (\boy - \byhat), & \label{M2-d budget} \\
&\sum_{j=1}^{m} (z^-_j+z^+_j) \le \Gamma, &  \label{M2-gamma} \\ 
& \bzero \le \bzn , \bzp \le \bone, & \label{M2-f z limits}
\\ & \bx,\by \ge \bzero. &
\label{M2-e budget}
\end{align}
\end{subequations}

The inner maximization in the objective function finds the worst-case cost and can be reformulated using Duality Theorems as shown in \cite{bertsimas2004price}. 
Parameter $\Gamma$ can take any value within $[0,m]$ where $m$ is the number of 
resource parameters that are subject to uncertainty. In particular, when $\Gamma =0$, problem~$\MB$ is equivalent to the nominal problem $\M$ with no uncertainty. As $\Gamma$ increases, it is expected for the solution to become more conservative against uncertainty. 
For a full budget ($\Gamma=m$), since $\mathcal{P}=\{i : c_{2,(i)} \ge \max_k \hspace{0.1cm} c_{1,(k)} \}$ is nonempty, as mentioned in Assumption \ref{assum:c}, the worst-case cost occurs when $\bzn=\bzero$ and $\bzp=\bone$.
Formulation $\MBm$ corresponds to model $\MB$ with full budget of uncertainty, which is equivalent to the traditional robust model where all constraints assume their worst-case 
values, independently.
\begin{subequations}\label{soy M2 with Budget}
\begin{align}
\MBm: \quad \min_{\bx,\by} \quad & \mathbf{c}_1'\bx+ \mathbf{c}_2'\big(
\boy-\by\big) , & \label{soy M2-a budget} \\
\text{s.t} \quad & \bA\bx+\bB\by \le \bg, & \label{soy M2-b budget}\\ & \by \le \boy, & \label{soy M2-d budget} \\ & \bx,\by \ge \bzero. &
\label{soy M2-e budget}
\end{align}
\end{subequations}

In what follows, we motivate our proposed ``effective budget-of-uncertainty'' approach. 
To this end, in Proposition~\ref{initial prop}, we first show that there exist cases in which after a certain threshold of $\Gamma$, increasing the budget of uncertainty would not affect the optimal solution of formulation $\MB$.  In the rest of this paper, we use $\bA_i$ and $\bB_i$ to denote the $i^{\text{th}}$ row of matrices $\bA$ and $\bB$, respectively, and use $g_i$ to show the $i^{\text{th}}$ element of the vector $\bg$.

\begin{proposition}\label{initial prop}
For budgets of uncertainty $\Gammaz \geq \Gammal$, problems $\MBz$ and $\MBl$  with optimal solutions $(\bxz, \byz,\bzzn, \bzzp)$ and $(\bxl, \byl, \bznl,\bzpl)$, respectively, result in the same optimal values $\bxz=\bxl$ and $\byz=\byl$ { if  
$\, \exists \, i \in \{1, \dots, n\} $ such that $c_{2,(i)} \ge \max_k \hspace{0.1cm} c_{1,(k)}$,  $\bA_i \bxl \geq 0$, and $\bA_i \bxl + \bB_i \byl = g_i$. }   
\end{proposition}
\proof{Proof.}
On one hand, due to the direction of the inner maximization in the objective function, the worst-case penalty cost occurs when there are only positive deviations from the nominal value (i.e. $\bzzn=\bznl=\bzero$). Therefore, if $\Gammaz \geq \Gammal$, 
 it must be that $\bzzp \geq \bzpl$ and 
$\byz \geq \byl$. 
%
%
{ On the other hand, 
since the $i^{\text{th}}$ constraint of~\eqref{M2-b budget} is binding for $\bxl$ and $\byl$ and $c_{2,(i)} \ge \max_k \hspace{0.1cm} c_{1,(k)}$, it means that it was not possible to further reduce the term $\bA_i\bx$ at optimality.} Hence, for any $\Gammaz \geq \Gammal$, we must have $\bA_i\bxz \geq \bA_i\bxl$. Therefore, we have $\bB_i\byz \leq g_i - \bA_i\bxz  \leq g_i - \bA_i \bxl$. Since $\bB \geq \bzero$, it can be concluded that $\byz \leq \byl$.
%
Hence, we conclude $\byz = \byl$, which based on the direction of objective function, also implies $\bxz=\bxl$ even though $\bzzp \ge \bzpl$. \Halmos
\endproof \vspace{1em}

Proposition \ref{initial prop} shows that it is possible that problem $\MB$ would become insensitive to increasing the budget of uncertainty after a certain threshold 
since constraints \eqref{M2-d budget} would become redundant and constraint(s)~\eqref{M2-b budget} would be binding. The intuition behind this idea is that due to constraint~\eqref{M2-b budget}, an increase in $\by$ must accompany a decrease in $\bA \bx$ in order for this constraint to remain feasible. If there is a lower bound on how much some elements of $\bA\bx$ can be decreased, then we would not be able to increase $\by$ by more than a certain level. Therefore, the available budget may not be entirely utilized as the budget increases. 
In particular, in power dispatch problems, there always exists a lower bound for the output of conventional generators (i.e., $\bA\bx$) due to technical issues. As the available wind power increases, the extra amount of available wind cannot be absorbed and must be curtailed. The absolute worst-case scenario occurs when all wind farms see the highest amount of wind within the uncertainty set, but there is often no change in optimal solution of the system to deal with the extra wind, other than paying the curtailment cost.
This observation, especially in practical SCED systems, is the mathematical motivation for the proposed approach. 

In what follows, we present a two-stage methodology for removing ineffective budgets of uncertainty and identifying the effective worst-case scenario which corresponds to budgets up to which the robust solution can be adjusted to address uncertainty. Consideration of any budget above this threshold will not affect the optimal robust solution. Using the proposed approach would help decision makers more intuitively study the trade-off between level on uncertainty and solution conservatism and would provide a less conservative approach for dealing with uncertainty.
We note that the proposed two-stage approach is not to be confused with the literature of two-stage robust optimization where some variables are determined in the first stage and some are decided in the second stage after the uncertain parameters are realized. Our two-stage approach consists of (1)~identifying subsets of the uncertainty set that include ineffective budgets of uncertainty, and (2)~formulating a robust problem for incorporating effective budgets and protecting against the effective worst-case only. 

The rest of this section is organized as follows. In Section~\ref{ch2 effective and admissible sets}, we first introduce admissible and effective 
intervals by taking the feasibility constraints into account. 
Then in Section~\ref{ch2 effective budget}, we identify an effective budget of uncertainty which will be used for finding a less conservative robust solution. Finally, in Section~\ref{two-stage model}, we summarize our two-stage robust optimization approach.


\subsection{Admissible Interval and Effective Worst-Case Scenario}\label{ch2 effective and admissible sets}
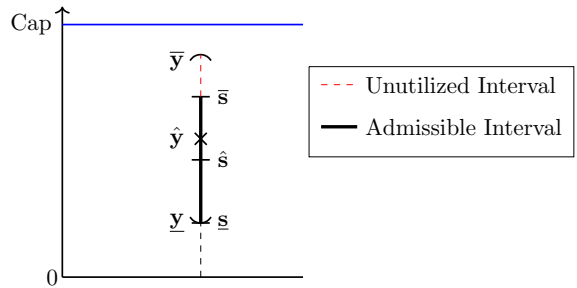
\begin{figure}[t]
\centering
\scalebox{0.8}
{\begin{circuitikz} 
\draw [thick, ->] (-1,1) -- (-1,5.5);
\draw [thick] (-1,1) -- (3,1);
\draw (-1.4,1) node [right] {0};
\draw (-2,5.2) node [right] {Cap};
\draw [thick, blue] (-1,5.2) -- (3,5.2);
\draw [dashed,red] (3.3,4.2) -- (3.9,4.2);
\draw (3.9,4.2) node [right] {Unutilized Interval};
\draw [ultra thick] (3.3,3.5) -- (3.9,3.5);
\draw (3.9,3.5) node [right] {Admissible Interval};
\draw (3.1,4.5) -- (7.5,4.5)--(7.5,3)--(3.1,3)--(3.1,4.5);
\draw [thick] (1.47,4.6) arc (30:150:0.2cm);
\draw (0.65,4.6) node [right] {$\boy$};
\draw [thick] (1.47,2) arc (-30:-150:0.2cm);
\draw (0.65,1.9) node [right] {$\buy$};
\draw [thick] (1.2,3.4) -- (1.4,3.2);
\draw [thick] (1.4,3.4) -- (1.2,3.2);
\draw (0.65,3.32) node [right] {$\byhat$};
\draw [thick] (1.15,4) -- (1.45,4);
\draw (1.45,4) node [right] {$\bos$};
\draw [thick] (1.15,1.9) -- (1.45,1.9);
\draw (1.45,1.9) node [right] {$\bus$};
\draw [thick] (1.15,2.95) -- (1.45,2.95);
\draw (1.45,2.965) node [right] {$\bshat$};
\draw [dashed](1.3,1) -- (1.3,4);
\draw [dashed, red] (1.3,4) -- (1.3,4.7);
\draw [ultra thick] (1.3,4) -- (1.3,1.9);
\end{circuitikz}}
\caption{\label{feasible window vs. uncertainty set} Comparison between the initial interval $[\buy,\boy]$, the admissible interval $[\bus,\bos]$.}
\end{figure}

We define the ``admissible'' interval $[\bus,\bos]$ as the largest subset of $[\bzero,\boy]$ such that for any value $\by \in [\bus,\bos]$
problem $\MB$ with a full budget of uncertainty ($\Gamma=m$) is always feasible \citep{wu2014robust}. 
We note that with a full budget, each element of variable $\by$ can take any value within its potential domain $[\bzero,\boy]$ independently. Figure~\ref{feasible window vs. uncertainty set} shows a schematic for the admissible interval. 
The ``effective worst-case scenario'' is the worst-case of the admissible interval for which the robust solution can be adjusted.
In what follows, we first present an optimization problem for finding the admissible interval, and then we identify the effective worst-case scenario, accordingly.

\subsubsection{Admissible 
Interval: }
When there is a full budget of uncertainty ($\Gamma =m$), each element of vector $\by$ in $\MBm$ can take any value within its potential domain in $[\bzero,\boy]$. 
The admissible interval $[\bus,\bos]$ is a subset of $[\bzero,\boy]$ (the 
lower and upper bounds for variable $\by$ in the full budget model) and may partially or entirely lay outside of the 
bounds of the uncertain parameters, i.e., $[\buy,\boy]$. 
Thus, for the upper and lower limits of the admissible interval, we have \citep{wu2014robust}:
\begin{subequations} \label{e8 relationships of Wmax Wmin and Smax Smin}
\begin{align}
\bos \le \boy, & \label{e8-a} \\ 
\bus \le \buy, & \label{e8-b}
\end{align}
\end{subequations}
where $\bos=\boy$ means that  interval $[\buy,\boy]$ is entirely admissible. 
Proposition \ref{prop first} demonstrates how to identify whether a given interval $[\bus,\bos]$ is entirely admissible. Remark \ref{remark 1} further explains how to find the largest such interval that is as close as possible to the limits of the interval $[\buy,\boy]$.

\begin{proposition} \label{prop first}
Any solution $\by \in [\bus, \bos]$ 
is a feasible solution 
of formulation $\MBm$ 
if there exists $\bx$, $\bus$, $\bos$, and $\balpha_i$ that meet the following conditions: 
\begin{subequations}\label{M3 upper bound}
    \begin{align} 
&\bA_i\bx+ \bB_i\bus +\bone'\balpha_i' \le g_i, & \forall i \in \{1, \dots, n \} \\
&\balpha_i' \ge \bB_i' \odot (\bos-\bus), & \forall i \in \{1, \dots, n \} \\ 
&\bx , \balpha_i' \ge \bzero. & \forall i \in \{1, \dots, n \}
    \end{align}
\end{subequations}
%
where $\balpha_i$ is the $i^{\text{th}}$ row of the matrix $\balpha$. 
\end{proposition}
The proof of Proposition \ref{prop first} is based on \cite{wu2014robust} and is provided in Appendix \ref{Appendix proof}.

\begin{remark} \label{remark 1}
The following optimization problem finds the largest 
admissible interval $[\bus,\bos]$ that has the smallest distance from the bounds of the interval $[\buy, \boy]$ of the uncertainty set $\cY$. 
\begin{subequations}\label{sub}
\begin{align} 
\min_{\bus,\bos, \balpha_i', \bx} \quad &  \bm{1}'(\boy-\bos)+\bm{1}'(\buy-\bus),  & \label{sub-a}\\
\mathrm{s.t.} \quad & \bA_i\bx+ \bB_i\bus +\bone'\balpha_i'   \le g_i, & \forall i \in \{1, \dots, n \} \label{sub-b}\\
&\balpha_i' \ge \bB_i' \odot (\bos-\bus), & \forall i \in \{1, \dots, n \} \label{sub-c}\\
& \bos \le \boy, & \label{sub-d}\\
&\bus \le \buy, &\label{sub-e}\\
&\bus \le \bos, &\label{sub-eee}\\
&\bx,\bus,\bos, \balpha_i' \ge \bzero. & \forall i \in \{1, \dots, n \} \label{sub-f} 
\end{align}
\end{subequations}
\end{remark}
In Remark~\ref{remark 1}, constraints \eqref{sub-b}, \eqref{sub-c}, and $\balpha_i \ge \bzero$ in \eqref{sub-f} identify 
whether a given interval $[\bus,\bos]$ is admissible. Constraints \eqref{sub-d} and \eqref{sub-e} demonstrate that the admissible interval may lay outside of the uncertainty set since the RHS parameters might not be fully utilized. Constraints \eqref{sub-eee} ensure that the admissible interval has an upper limit greater or equal to its lower limit. Given that the objective function $\eqref{sub-a}$ minimizes the gap between $[\buy,\boy]$ and $[\bus,\bos]$, 
formulation \eqref{sub} always obtains the largest possible admissible interval that has the smallest distance from the initial uncertainty set. In what follows, we describe important characteristics of the admissible interval. 
\begin{proposition} \label{prop - admissible interval}
For any $j$, the admissible interval $[\suj,\soj]$ 
can always be categorized as one of the following four cases:
\begin{enumerate}
\item $\suj=\yuj$ and \hspace{0.05cm} $\soj=\yoj$
\item $\suj=\yuj$ and \hspace{0.05cm} $\yhj \le \soj < \yoj$
\item $\suj=\yuj$ and \hspace{0.05cm} $\yuj < \soj < \yhj$
\item $\suj=\soj \le \yuj$
\end{enumerate}
\end{proposition}
\proof{Proof.}
Let $\alpha_{ij}$, $\suj$, and $\soj$ denote the $j^{\text{th}}$ element of $\balpha_i$, $\bus$, and $\bos$, respectively, and $B_{ij}$ denote the element of row $i$ and column $j$ in matrix $\bB$. From \eqref{sub-c}, $\alpha'_{ij} \ge B_{ij} (\soj-\suj),~\forall i,j$, if $\alpha_{ij}'$ is feasible. By Assumption \ref{assum:B}, $B_{ij} \ge 0$. 
Let us separate the case of $B_{ij}=0$ and $B_{ij}>0$: 

\noindent
(i) If $B_{ij}=0$, $\alpha'_{ij}  \ge B_{ij} (\soj-\suj)$ would become redundant $\forall \soj, \suj$. 
Thus, due to the minimization objective function, $\soj=\yoj$ and $\suj=\yuj$ at optimality which corresponds to case (a). 
%

\noindent
(ii) If, on the other hand, $B_{ij}>0$, we again separate the cases in which $\alpha'_{ij}>0$ or $\alpha'_{ij}=0$. 

\noindent 
First, for a particular $i$ and $j$ corresponding to a feasible $\alpha'_{ij}>0$, we can write \eqref{sub-b} as %
\begin{equation}\label{exp}
\bA_i\bx + \sum_{\substack{k=1 \\ k\neq j}}^{m} (B_{ik}\underline{s}_k +\alpha'_{ik}) + B_{ij}\suj +\alpha'_{ij} \le g_i\,.
\end{equation}%
Since $\alpha'_{ij}>0$, constraint~\eqref{sub-c} implies $\suj < \soj$. Substituting $\alpha'_{ij}$ from \eqref{sub-c} in \eqref{sub-b}, we have $B_{ij}\suj + B_{ij}(\soj-\suj)=B_{ij}\soj \le B_{ij}\suj + \alpha'_{ij}$. Thus, inequality \eqref{exp} holds $\forall \suj \in [0,\yuj]$, since $B_{ij}\suj < B_{ij}\soj \le B_{ij}\suj + \alpha'_{ij}$, $\forall \alpha'_{ij}>0$.  
Due to the minimization objective function, $\suj=\yuj$ at optimality. Depending on the value of a feasible $\alpha'_{ij}$, by substitution we observe from \eqref{sub-b} and \eqref{sub-c} that (note $\suj=\yuj$): 
\begin{itemize} \itemsep 1pt
\item If $B_{ij}(\yoj-\yuj) \le \alpha'_{ij}$, then $B_{ij}\suj + B_{ij}(\yoj-\yuj)=B_{ij}\yoj \le B_{ij}\suj + \alpha'_{ij}$. The later inequality holds $\forall \soj \in [0,\yoj]$ since $B_{ij}\soj \le B_{ij}\yoj$ . Thus, 
$\soj=\yoj$ at optimality, which corresponds to case (a). 
\item Similarly, by substitution we can show, if $B_{ij}(\yhj-\yuj) \le \alpha'_{ij} < B_{ij}(\yoj-\yuj)$, then $\yhj \le \soj < \yoj$ at optimality, which corresponds to case (b).
\item If $0 < \alpha'_{ij} < B_{ij}(\yhj-\yuj)$, then $\yuj < \soj < \yhj$ at optimality, corresponding to case (c).
\end{itemize}

\noindent
Second, for $\alpha'_{ij}=0$, 
we show $\suj=\soj \le \yuj$ by contradiction: Let sets $M=[0,\yuj]$ and $N=(\yuj,\yoj]$ where $M \cap N = \emptyset$
and $M \cup N = [0, \yoj]$. Assume $\exists\,\, \suj,\soj \in N$, such that $\suj=\soj$. This contradicts constraint $\suj\le\yuj$. Thus 
$\suj,\soj \in M$ meaning that $\suj=\soj \le \yuj$, which corresponds to case (d).
\Halmos \endproof  \vspace{1em}

So far, we identified the admissible interval by removing 
the part of interval $[\buy,\boy]$ that can never be utilized. Now, consider the following uncertainty set $\cU^A$ with a budget $\Gamma^A$ based on the admissible interval as follows:
\begin{equation}
\label{master problem uncertainty set}
\cU^{A} = \Big\{ \hspace{0.051cm} \bts \in \mathbb{R}^{{m}} : \hspace{0.2cm}  \bts = \bshat+\bz^+\odot(\bos-\bshat)+\bz^-\odot(\bus-\bshat), \hspace{0.2cm} \sum_{j=1}^{m} (z^-_{j}+z^+_{j}) \le \Gamma^A, \hspace{0.2cm} \bzero \le\bz^+,\bz^- \le \bone \hspace{0.1cm} \Big\},
\end{equation}
where $\bshat$ is the middle point of the interval $[\bus,\bos]$, and $\bz^+$ and $\bz^-$ denote positive and negative scaled deviations from $\bshat$, respectively. Consider formulation $\MAA$ where the uncertainty set $\cU$ in formulation $\MB$ is substituted with the uncertainty set $\cU^A$.
\begin{subequations}\label{admissible robust}
\begin{align}
\MAA: \quad \min_{\bx,\by} \quad & \Big\{ \mathbf{c}_1'\bx+ \max_{\bz} \mathbf{c}_2'\Big(
\bshat + \bzn \odot (\bus - \bshat) + \bzp \odot (\bos - \bshat)-\by\Big) \Big\}, & \label{M2-a budget A} \\
\text{s.t} \quad & \bA\bx+\bB\by \le \bg, & \label{M2-b budget A}\\ & \by \le \bshat + \bzn \odot (\bus - \bshat) + \bzp \odot (\bos - \bshat), & \label{M2-d budget A} \\
&\sum_{j=1}^{m} (z^-_j+z^+_j) \le \Gamma^A, &  \label{M2-gamma A} \\ 
& \bzero \le \bzn , \bzp \le \bone, & \label{M2-f z limits A}
\\ & \bx,\by \ge \bzero. &
\label{M2-e budget A}
\end{align}
\end{subequations}

When there is a full budget, due to the direction of the objective function, $\bzn=0$ and $\bzp=1$. Thus, the following formulation corresponds to $\MAm$ with a full budget $\Gamma=m$.
\begin{subequations}\label{ad soy M2 with Budget}
\begin{align}
\MAm: \quad \min_{\bx,\by} \quad & \mathbf{c}_1'\bx+ \mathbf{c}_2'\big(
\bos-\by\big) , & \label{ad soy M2-a budget} \\
\text{s.t} \quad & \bA\bx+\bB\by \le \bg, & \label{ad soy M2-b budget}\\ & \by \le \bos, & \label{ad soy M2-d budget} \\ & \bx,\by \ge \bzero. &
\label{ad soy M2-e budget}
\end{align}
\end{subequations}
The following Remark shows that 
problems $\MBm$ and $\MAm$ have the same feasible values for $(\bx,\by)$, which means that the admissible interval would not remove any solution $(\bx,\by)$ from the feasible region of problem $\MBm$.
\begin{remark} \label{remark RB equals RA}
Problems $\MBm$ and $\MAm$ have the same feasible space for $\bx$ and $\by$.
\end{remark} %
\proof{Proof.}
Let $(\bx,\by, \bzn, \bzp) \in \mathbf{X}$ and $(\bx^A,\by^A, \bznA, \bzpA) \in \mathbf{X}^A$ be feasible solutions of $\MBm$ and $\MAm$, respectively. Since there is a full budget in both $\MBm$ and $\MAm$, $\bzn=\bznA=\bzero$ and $\bzp=\bzpA=\bone$ due to the direction of the objective function. To conclude that the two formulations correspond to the same feasible space for $\bx$ and $\by$, it is sufficient to show that any feasible solution $(\bx,\by,\bzero,\bone)$ in $\mathbf{X}$ is also feasible for $\mathbf{X}^A$ and vice-versa.

\noindent
First, let $(\bx^A,\by^A, \bzero, \bone) \in \mathbf{X}^A$, where $\by^A \le \bos$. Since $\bos$ falls into one of the cases of Proposition \ref{prop - admissible interval}, it satisfies the conditions of Proposition \ref{prop first}, meaning that $\bos \le \boy$. Thus, $(\bx^A,\by^A,\bzero,\bone) \in \mathbf{X}$ as well.
\noindent
Similarly, let $(\bx,\by,\bzero,\bone) \in \mathbf{X}$, where $\by \le \boy$. We observe that for some values of $\boy$ which cannot be utilized, constraint \eqref{M2-d budget} becomes redundant.
Thus, $\by \le \bos$ and  $(\bx,\by,\bzero,\bone) \in \mathbf{X}^A$.
\Halmos \endproof  \vspace{1em}

From Remark \ref{remark RB equals RA}, we can conclude that for $\Gamma = \Gamma^A = m$, both formulation $\MB$ and $\MAA$ have the same worst-case performance and correspond to the same solutions $\bx$ and $\by$ for the full budget of uncertainty. 
Thus, to find a less conservative robust solutions, we can use the effective worst-case scenario, which falls within the admissible interval, as opposed to the absolute worst-case scenario. 

\subsubsection{Effective 
Worst-Case Scenario:}
We refer to the worst-case scenario of admissible interval as the effective worst-case scenario and identify a subset of $[\bus,\bos]$ within which the effective worst-case would always occur.
\begin{proposition} \label{prop worst-case realization interval}
The worst-case realization of the admissible uncertainty set ${\cU^{A}}$ always occurs within interval $[\bshat,\bos]$ where $\bshat$ is the middle point of the interval $[\bus,\bos]$.
\end{proposition}
\proof{Proof.} 
Regardless of the value of the budget parameter, using the definition of $\cU^{A}$, we note that $[\bus,\bos]=[\bus,\bshat] \cup [\bshat,\bos]$. Thus, we can re-write the worst-case of constraint \eqref{e5}, presented in Appendix \ref{Appendix proof}, by considering
$\bm{\beta}(\bz^+,\bz^-)$  instead of $\bm{\beta}_i(\br)$, 
where: %
\begin{subequations} \label{inner max stilde}
\begin{align}
\bm{\beta}_i(\bz^+,\bz^-)  = \max_{\bz^+,\bz^-} \quad &
\bB_i \Big( \bz^+\odot(\bos-\bshat)+\bz^-\odot(\bus-\bshat)\Big), & \forall i \in \{1, \dots, n\} \label{stilde_proof-1}\\ 
\text{s.t. } \quad & \bzero \le \bz^+ \le \bone, & \label{stilde_proof-2 +} \\
& \bzero \le \bz^- \le \bone. & \label{stilde_proof-2 -}
\end{align}
\end{subequations}
Let $\bm{\eta}$ and $\bm{\zeta}$ be matrices with their columns $\bm{\eta}_i'$ and $\bm{\zeta}_i'$ being dual variables corresponding to constraints \eqref{stilde_proof-2 +} and \eqref{stilde_proof-2 -}, respectively, for each row $i$. 
The dual formulation of \eqref{inner max stilde} is: 
\begin{subequations} \label{pe6}
\begin{align}
\min_{\bm{\eta}_i',\bm{\zeta}_i'}\quad & \bone'\bm{\eta}_i' + \bone'\bm{\zeta}_i', & \forall i \in \{1, \dots, n \}  \label{pe6-1} \\
\text{s.t. } \quad &\bm{\eta}_i' \ge \bB_i'\odot(\bos-\bshat ),  & \forall i \in \{1, \dots, n \}  \label{pe6-2} \\
&\bm{\zeta}_i' \ge \bB_i'\odot(\bus-\bshat), & \forall i \in \{1, \dots, n \} \label{pe6-3} \\ 
& \bm{\eta}_i', \bm{\zeta}_i'\ge \bzero.& \forall i \in \{1, \dots, n \}   
\end{align}
\end{subequations}
Since $\bus-\bshat \le \bzero$ and $\bB_i \geq \bzero$, 
constraint \eqref{pe6-3} is redundant and thus $\bm{\zeta}_i'=\bm{0}$. 
We remove the redundant constraint \eqref{pe6-3} and take the dual of model \eqref{pe6} again for each $i$ where vector $\mathbf{p}_i'$ is a column of dual variables corresponding to constraint \eqref{pe6-2}. Hence, we obtain the following formulation: 
\begin{subequations} \label{inner max stilde the second}
\begin{align}
\max_{\mathbf{p}_i'} \quad &
\bB_i \Big( \mathbf{p}_i' \odot(\bos-\bshat)\Big), & \forall i \in \{1, \dots, n \} \\ 
\text{s.t. } \quad & \bzero \le \mathbf{p}_i' \le \bone, & \forall i \in \{1, \dots, n \}
\end{align}
\end{subequations}
which is  equivalent to \eqref{e5-2} and shows that $\bts \in [\bshat,\bos]$.
\Halmos \endproof  \vspace{1em}
%
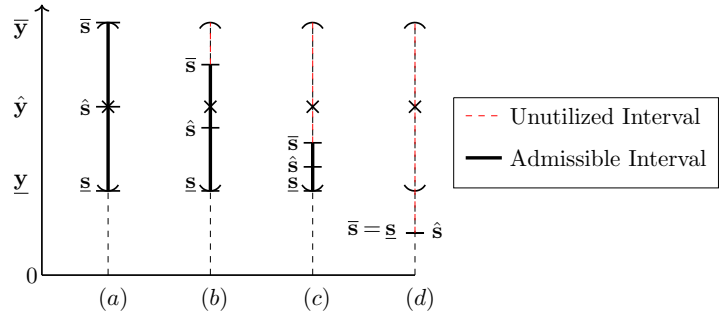
\begin{figure}[t]
\centering
\scalebox{0.8}
{\begin{circuitikz} 
\draw [thick, ->] (-1.5,0.5) -- (-1.5,5);
\draw [thick] (-1.5,0.5) -- (4.7,0.5);
\draw (-1.9,0.5) node [right] {0};
\draw [thick] (-0.23,4.6) arc (30:150:0.2cm);
\draw [thick] (1.47,4.6) arc (30:150:0.2cm);
\draw [thick] (3.17,4.6) arc (30:150:0.2cm);
\draw [thick] (4.87,4.6) arc (30:150:0.2cm);
\draw (-2.1,4.6) node [right] {$\boy$};
\draw [thick] (-0.23,2) arc (-30:-150:0.2cm);
\draw [thick] (1.47,2) arc (-30:-150:0.2cm);
\draw [thick] (3.17,2) arc (-30:-150:0.2cm);
\draw [thick] (4.87,2) arc (-30:-150:0.2cm);
\draw (-2.1,2) node [right] {$\buy$};
\draw [thick] (-0.5,3.2) -- (-0.3,3.4);
\draw [thick] (-0.5,3.4) -- (-0.3,3.2);
\draw [thick] (1.2,3.2) -- (1.4,3.4);
\draw [thick] (1.2,3.4) -- (1.4,3.2);
\draw [thick] (2.9,3.2) -- (3.1,3.4);
\draw [thick] (2.9,3.4) -- (3.1,3.2);
\draw [thick] (4.6,3.2) -- (4.8,3.4);
\draw [thick] (4.6,3.4) -- (4.8,3.2);
\draw (-2.1,3.3) node [right] {$\byhat$};
\draw [thick] (-0.6,4.7) -- (-0.2,4.7);
\draw [thick] (1.15,4) -- (1.45,4);
\draw [thick] (2.85,2.7) -- (3.15,2.7);
\draw [thick] (4.55,1.2) -- (4.85,1.2);
\draw (-1,4.65) node [right] {$\bos$};
\draw (0.75,4) node [right] {$\bos$};
\draw (2.45,2.7) node [right] {$\bos$};
\draw (3.45,1.25) node [right] {$\bos=\bus$};
\draw [thick] (-0.6,1.9) -- (-0.2,1.9);
\draw [thick] (1.15,1.9) -- (1.45,1.9);
\draw [thick] (2.85,1.9) -- (3.15,1.9);
\draw [thick] (4.55,1.2) -- (4.85,1.2);
\draw (-1,2) node [right] {$\bus$};
\draw (0.75,2) node [right] {$\bus$};
\draw (2.45,2) node [right]
{$\bus$};
\draw [thick] (-0.6,3.3) -- (-0.2,3.3);
\draw [thick] (1.15,2.95) -- (1.45,2.95);
\draw [thick] (2.85,2.3) -- (3.15,2.3);
\draw [thick] (4.55,1.2) -- (4.85,1.2);
\draw (-1,3.3) node [right] {$\bshat$};
\draw (0.75,2.95) node [right] {$\bshat$};
\draw (2.45,2.35) node [right]
{$\bshat$};
\draw (4.85,1.25) node [right]
{$\bshat$};
\draw [dashed](-0.4,0.5) -- (-0.4,4.7);
\draw [dashed](1.3,0.5) -- (1.3,4.7);
\draw [dashed](3,0.5) -- (3,4.7);
\draw [dashed](4.7,0.5) -- (4.7,4.7);
\draw [dashed, red] (1.3,4) -- (1.3,4.7);
\draw [dashed, red] (3,2.7) -- (3,4.7);
\draw [dashed, red] (4.7,1.2) -- (4.7,4.7);
\draw [ultra thick] (3,2.7) -- (3,1.9);
\draw [ultra thick] (1.3,4) -- (1.3,1.9);
\draw [ultra thick] (-0.4,4.7) -- (-0.4,1.9);
\draw (-0.7,0.1) node [right]
{($\boldmath{a}$)};
\draw (1,0.1) node [right]
{($\boldmath{b}$)};
\draw (2.7,0.1) node [right]
{($\boldmath{c}$)};
\draw (4.4,0.1) node [right]
{($\boldmath{d}$)};
\draw (5.2,2.7) node [right]{\begin{circuitikz}
\draw [dashed,red] (3.3,4.2) -- (3.9,4.2);
\draw (3.9,4.2) node [right] {Unutilized Interval};
\draw [ultra thick] (3.3,3.5) -- (3.9,3.5);
\draw (3.9,3.5) node [right] {Admissible Interval};
\draw (3.1,4.5) -- (7.5,4.5)--(7.5,3)--(3.1,3)--(3.1,4.5);
\end{circuitikz}};
\end{circuitikz}}
\caption{Comparison of the initial interval $[\buy,\boy]$ and the admissible interval $[\bus,\bos]$ for the four possible cases of Proposition~\ref{prop - admissible interval}.} 
\label{feasible window_new}
\end{figure}
Proposition \ref{prop worst-case realization interval} shows the equivalency of models with $[\bus,\bos]$ and $[\bshat,\bos]$ as uncertain intervals in terms of their worst-case scenarios. Since the effective worst-case scenario only corresponds to positive scaled deviations from $\bshat$, we can conclude $\bzn=\bzero$ in \eqref{master problem uncertainty set}. We will use this result in next section to develop our two-stage approach.

\noindent
In what follows, we identify an \emph{effective budget of uncertainty} ($\Gamma^E$) to be incorporated in the problem using the definitions made so far. 

\subsection{Effective Budget of Uncertainty} \label{ch2 effective budget}
In the conventional budget approach, the budget $\Gamma$ 
controls the sum of scaled deviations from the nominal value $\byhat$, which is the middle point of the 
interval $[\buy, \boy]$.  In Proposition~\ref{prop worst-case realization interval}, we showed that the effective worst-case scenario always occurs within the interval $[\bshat,\bos]$, i.e., $\bzn=\bzero$, and that we can replace interval $[\buy,\boy]$ with $[\bshat,\bos]$ 
 without changing the worst-case performance for full budget. 
 However, when there is no full budget, 
 the deviations in $[\bshat,\bos]$ are captured from the nominal value 
$\bshat$ which is often not equal to $\byhat$ as shown in Fig \ref{feasible window_new}. 
So, the conventional definition of $\Gamma$ may need to be adjusted.
In this section, we propose a new definition of effective budget of uncertainty on 
the set $\cU^A$ with the following properties:
\begin{enumerate}
    \item For $\Gamma=0$ and $\Gamma=m$ (zero budget and full budget, respectively), it generates the same solutions as those of the conventional budget approach.
    \item For other values of $0<\Gamma<m$, depending on the corresponding case from Proposition~\ref{prop - admissible interval}, it finds the effective worst-case scenario within the budget $\Gamma$. 
\end{enumerate}

\begin{figure}[t] 
    \centering
  \begin{tikzpicture}[font=\small]
\pgfplotsset{width=7cm, height=5.4cm}
  \begin{axis}[legend image post style={scale=0.6},legend style={{font=\fontsize{8.7}{8.7}\selectfont}, at={(1.04,0.75)}, anchor=north west},
   xmin=-0.05,
   xmax=1.05,
   xtick={0,1},
   xlabel= {Scaled deviation},
   ymin=0,
   ymax=50,
    ytick={0,10,30,50},
    yticklabels={$\shj$,$\yhj$,$\soj$,$\yoj$},
    ylabel={Uncertain parameter}
  ] 
\addplot [no marks, gray!15,name path=yhat] plot coordinates {
(0,10) (1,10)};
\addplot [no marks, gray!15,name path=sbar] plot coordinates {
(0,30) (1,30)};
\addplot [no marks, black,name path=ybar] plot coordinates {
(0,50) (1,50)};
\addplot [no marks, black,name path=zero] plot coordinates {
(0,0) (1,0)};
\addplot[gray!15] fill between[ 
    of = yhat and zero];
\addplot[gray!15] fill between[ 
    of = ybar and sbar];
\addplot [no marks, black,dashed,thick, name path=conventional approach] plot coordinates {
(0,10)(1,50)};
%
\addplot [no marks, black,thick, name path=effective] plot coordinates {
(0,0)(1,30)};
\addplot [no marks, green!80!black,thick, name path=proposed approach] plot coordinates {
(0,10)(1,30)};
\addplot[scatter/classes={
a={mark=*,black},
b={mark=*,green!80!black}},
scatter,only marks,
scatter src=explicit symbolic] coordinates {(0,0) [a] (0,10) [b] (1,50) [a] (0.5,30) [a] (0.33,10) [a] (1,30) [b]};
 \node[green!80!black] at (axis cs:0, 13.3) {$A$};
\node at (axis cs:0, 3.3) {$A'$};
\node at (axis cs:0.5, 33.3) {$C$};
\node at (axis cs:0.33, 13.3) {$C'$};
\node at (axis cs:1, 46.5) {$B$};
\node[green!80!black] at (axis cs:1, 26.5) {$B'$};
\legend{,,,,,Ineffective deviation, Original interval, Admissible interval, Effective Interval}
  \end{axis}
\end{tikzpicture}
\caption{Comparison of the uncertain parameters in case (b) for Original interval, Admissible interval, Effective Interval. Thu budget of uncertainty is 1 unit and is based on the scaled deviation of each uncertainty set.}
\label{Mapping Fig}
\end{figure}
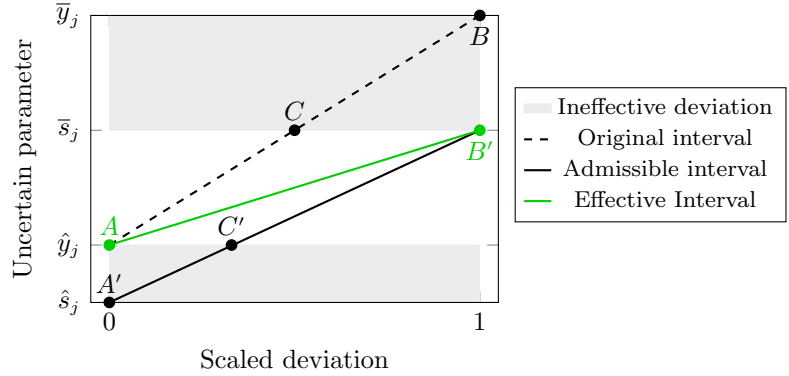

Before presenting the definition of an effective budget, let us first present an intuition behind what is expected. Figure~\ref{Mapping Fig} compares the uncertain parameters in the original, admissible, and effective intervals of case (b) based on their scaled deviations for 1 unit of the budget of uncertainty. The shaded regions in Figure~\ref{Mapping Fig} show the areas within which the budget of uncertainty of the conventional approach does not impact the solutions. Line $AB$ corresponds to the deviations of the uncertain parameter in the original interval $[\yhj,\yoj]$ using the conventional approach. However, only deviations corresponding to line segment $AC$ are effective since interval $(\soj,\yoj]$ cannot be utilized. Similarly, line $A'B'$ corresponds to the uncertain parameter of the admissible interval $[\shj,\soj]$. Recall that the conventional approach with uncertainty set $\cU$ and zero budget of uncertainty ($\Gamma=0$) allows variable $y_j$ to take values up to its nominal value $\yhj$ as long as it satisfies all other constraints. To preserve this point in the admissible interval $[\shj,\soj]$, we need to allow variable $y_j$ to take values up to $\yhj$ for a zero budget in $\cU^A$. Thus, deviations corresponding to line segment $A'C'$ are not effective. Finally, line $AB'$ corresponds to deviations that are both admissible and effective and are adjusted based on the length of admissible interval $[\shj,\soj]$. 
Thus, line $AB'$ corresponds to all robust solutions whose level of level of conservatism can be adjusted by changing the 
budget of uncertainty. 

For the other cases, a similar argument can be made. In cases (c) and (d), a zero budget of uncertainty ($\Gamma = 0$) in the conventional approach corresponds to solution $y_j = \soj$ since $\soj < \yhj$. To generate the same solution in $\cU^{A}$ for such cases, a one-directional positive deviation $z^+_j$ from the new nominal value $\shj$ is required since $y_j \le \shj+z^+_j(\soj-\shj)$. Otherwise, a zero budget in $\cU^{A}$ would result in $y_j=\shj$ (as $z^+_j=0$), which is different from the solution of the conventional approach with no uncertainty. Thus, to keep the properties of the conventional approach, we need to allow such deviations in the proposed approach but must not take them into account in the effective budget constraint for cases (c) and (d), since this amount of deviation has no effect on the solution.


Based on the properties discussed above, in Definition \ref{Def1}, we formally present an effective uncertainty set $\cU^{E}$ with an effective budget $\Gamma^E$. Finally, Theorem~\ref{prop:effective} 
summarizes the properties of the proposed effective budget of uncertainty. 
\begin{definition} \label{Def1}
$\cU^{E}$ is the effective uncertainty set with an effective budget of uncertainty $\Gamma^{E}$ and is defined as
\begin{equation}
\label{Uncertainty set - stage 2}
\cU^{E} = \Big\{ \hspace{0.2cm} \bts \in \mathbb{R}^{{m}} : \quad  \bts = {\bshat} + \br\odot({\bos}-{\bshat}), \quad \sum_{j=1}^{m} e_{j}r_{j} \le \Gamma^{E}, \quad \bv \le \br \le \bone. \hspace{0.2cm} \Big\} 
\end{equation}
where vectors $\bv$, $\bbee$ 
and the new budget $\Gamma^{E}$ are parameters calculated based on existing information as follows, and $\sign(a)$ is a sign function that returns 1 if parameter $a\ge 0$, and -1 otherwise. 
\begin{subequations} \label{define components of the S II uncertainty set}
\begin{align}
\displaystyle 
& \bh=\frac{1}{2}\Big(1-\sign(\byhat-\bos)\Big), \\[4pt]
& \displaystyle \bbee=\bh\odot (\frac{\bos-\bshat}{\boy-\byhat}), \\[10pt] & \displaystyle  
\bv=\bh\odot(\frac{\byhat-\bshat}{\bos-\bshat}), \\[10pt]
& \displaystyle  \Gamma^{E} = \Gamma + \sum_{j=1}^{m} v_{j} (\frac{\overline{s}_{j}-\hat{s}_{j}}{\overline{y}_{j}-\hat{y}_{j}}). 
\end{align}
\end{subequations}
\end{definition}
\begin{theorem}\label{prop:effective}
The effective uncertainty set $\cU^{E}$ only considers deviations within the effective interval $[\byhat,\bos]$ into the effective budget $\Gamma^{E}$.   
\end{theorem}
\proof{Proof.}
Considering the definition of $h_j$ and $v_j$, Table shows their values for all cases as follows.
\begin{table}[!h]
\centering
\footnotesize
\caption{Values of parameters $h_j$ and $v_j$ for all 4 cases}
\label{Table: parameters}
\begin{tabular}{lccccc}
\hline 
Parameters&& case (a)&case (b)&case (c)&case (d) \\ \hline
$h_j=\frac{1}{2}\big(1-\sign(\yhj-\soj)\big)$ && 1&1&0&0 \\
$v_j=h_j
(\yhj-\shj)/(\soj-\shj)$ && 0&non-zero&0&0
\\ \hline
\end{tabular}
\end{table}
Since $e_j$ is a function of $h_j$, the budget constraint in \eqref{Uncertainty set - stage 2} 
only considers effective deviations by letting the required deviations of cases (c) and (d) happen but not allocating an effective budget of uncertainty for such deviations.

\noindent
On the other hand, 
the nonzero value of $v_j$ is the scaled deviation required to map $\shj$ to $\yhj$ in case (b). Thus, constraint $v_j \le r_j$ ensures $\shj$ is mapped to $\yhj$. To ensure the proposed approach allows this mapping without using the budget of uncertainty, $\Gamma$ is linearly mapped to $\Gamma+\sum_{j=1}^{m}
v_j$, where 
\begin{equation} \label{e6 for proof}
\sum_{j=1}^{m}
h_jr_{j} \le \Gamma + \sum_{j=1}^{m}
v_j
\end{equation}
\noindent
Note that $v_j$ and $r_j$ are scaled deviations based on the length of $(\soj-\shj)$, while $\Gamma$ is a scaled parameter based on the magnitude of $(\yoj-\yhj)$. To normalize the scaled deviations, the factor $\frac{\soj-\shj}{\yoj-\yhj}$ is multiplied by $v_j$ and $r_j$. Doing so, formulations \eqref{Uncertainty set - stage 2} and \eqref{define components of the S II uncertainty set} can be recovered.
\Halmos \endproof \vspace{1em}

We note that the range of values that $\Gamma$ and $\Gamma^{E}$ can take are different. Recall that in the conventional budget approach, $\Gamma$  takes a value within $[0,m]$  \citep{bertsimas2004price}.
In the proposed approach, however, the effective budget 
$\Gamma^{E} \in \Big[ \sum_{j=1}^{m} v_{j} (\frac{\overline{s}_{j}-\hat{s}_{j}}{\overline{y}_{j}-\hat{y}_{j}}) \,~, \,~\sum_{j=1}^{m} v_{j} (\frac{\overline{s}_{j}-\hat{s}_{j}}{\overline{y}_{j}-\hat{y}_{j}})+m \Big]$. The term $\sum_{j=1}^{m} v_{j} (\frac{\overline{s}_{j}-\hat{s}_{j}}{\overline{y}_{j}-\hat{y}_{j}}) \ge 0$ is constant and is only used for case (b) as shown before.
Particularly, $\Gamma=0$ in the conventional approach would generate the same solution as $\Gamma^{E}=\sum_{j=1}^{m} v_{j} (\frac{\overline{s}_{j}-\hat{s}_{j}}{\overline{y}_{j}-\hat{y}_{j}})$ in the proposed approach
(point $A$ of Figure~\ref{Mapping Fig}). 
On the other hand, when the entire uncertainty set is admissible, 
(i.e., case (a)), considering the parameters of Table \ref{Table: parameters} 
we have $\Gamma^{E}=\Gamma$, and $\cU^E = \cU$. Thus, the proposed approach produces becomes equivalent to the conventional budget approach. 

\subsection{The Proposed Two-Stage Approach} \label{two-stage model} 
Based on the definition of $\cU^{E}$, we now propose a less conservative two-stage robust approach to solve the previously-described robust problems with 
resource uncertainty. By incorporating only the effective budget of uncertainty in the model, we optimize for the effective worst-case scenario instead of the absolute worst-case scenario of uncertainty. 

{ \bf Stage (I):} 
In this stage, we solve the auxiliary optimization problem \eqref{sub} 
to find the admissible interval 
$[\bus,\bos]$ that has the smallest distance from 
interval $[\buy,\boy]$.
We then consider the effective worst-case scenario to form $[\bshat,\bos]$ and calculate parameter $\Gamma^{E}$ and vectors $\bv$ and $\bbee$ using the set of equations in \eqref{define components of the S II uncertainty set}.

{\bf Stage (II):} In this stage, we use the output of Stage (I) and solve the following optimization problem~\eqref{master problem - stage II} which incorporates the effective budget of uncertainty into the robust model.
\begin{subequations}
\label{master problem - stage II}
\begin{align} 
\min_{\bx,\by} \quad & \Big\{ \mathbf{c}_1'\bx+ \max_{\br} \mathbf{c}_2'\Big(\bshat+\br \odot (\bos-\bshat)-\by \Big) \Big\}, & \\ \text{s.t} \quad & \bA\bx+ \bB\by \le \bm{g}, & \\ & \by \le {\bshat} + \br \odot ({\bos}-{\bshat}), & \\ & \sum_{j=1}^{m} e_{j}r_{j} \le \Gamma^{E}, & \\ & \bv \le \mathbf{r} \le \bone, &\\&
\bx,\by \ge \bzero. & 
\end{align}
\end{subequations}

\noindent
The proposed two-stage approach provides insights on the trade-off between robustness and solution conservatism in problems with resource 
uncertainty.
We note that to set the desired level of budget, one can start with the same budget $\Gamma$ as the \citet{bertsimas2004price} approach and then find the corresponding effective worst-case scenario. This would eliminate the ineffective budget while preserving the solution robustness.
From a managerial point of view, this approach 
allows the users to determine the level of budget of uncertainty that can impact the robust solution. 
In addition, decision makers can explicitly observe how increasing the budget of uncertainty would lead to a more conservative robust solution, as intuitively expected. This would allow decision makers to optimize the system under the effective worst-case scenario for any budget, as opposed to using the 
the traditional approach where the optimal solution may not change for higher budgets and the solution may be overly conservative.
{
\subsection{Extension to Other Uncertainty Sets}\label{sec:extension}
So far, we considered an interval uncertainty on $\bty$ in the form of upper bounds and lower bounds on each element of $\bty$, for simplicity. In this section, we demonstrate that the proposed approach can be applied to any shape of uncertainty, as long as the underlying set is convex, bounded, and closed. In this section, we demonstrate the overall approach for such a general uncertainty set. We later provide the specific steps for an example using an ellipsoidal uncertainty in Appendix~\ref{Appendix ellipsoid}.

Let $\cY$ be a nonempty bounded closed convex set in $\mathbb{R}^m$ that represents the uncertainty set on the parameter $\bty$. Since $\cY$ is bounded, for each element of the uncertain parameter $\tilde{y}_j$: 
\begin{equation}
  \exists \, \underline{y}_j, \, \overline{y}_j  \quad \text{such that} \quad \hspace{0.1cm} \yuj \le \, \tilde{y}_j \, \le \yoj, \qquad  \forall \tilde{y}_j \in \cY
    \end{equation} 
Let $\cH$ denote the axis-aligned minimum bounding box on the uncertainty set $\cY$ defined as follows:

\vspace{1em}
\begin{definition}
The set $\cH =\Big\{ \hspace{0.051cm} \bty \in \mathbb{R}^{{m}} : \hspace{0.2cm}  \buy \le \bty \le \boy \Big\}$ is the \textit{axis-aligned minimum bounding box} of a nonempty bounded closed convex set $\cY$ where $\forall j: \, \, \overline{y}_j = \underset{j}{\max}\{\tilde{y}_j \in \cY \}$ and $\underline{y}_j = \underset{j}{\min}\{\tilde{y}_j \in \cY\}$.
\end{definition} 
\vspace{1em}
The axis-aligned minimum bounding box $\cH$ provides the tightest upper and lower bounds for the uncertainty set $\cY$, regardless of the shape of $\cY$. Based on this definition, the following steps will be used in determining the admissible and effective intervals for any budget-based uncertainty set $\cU$ that is built on $\cY$: 

\begin{enumerate}
    \item Find $\cH$, the axis-aligned minimum bounding box of $\cY$. 
    \item Within $\cH$, find the admissible interval based on Remark~1 and Proposition~3. Then, intersect the admissible interval with $\cY$ to find the admissible uncertainty set, $\cU^A$.
    \item Find the effective part of the admissible interval depending on the shape of the uncertainty set for $\cU$, similar to Proposition \ref{prop worst-case realization interval}. 
\end{enumerate}
A hypothetical two-dimensional schematic of these steps for a general convex, bounded, and closed uncertainty set $\cY$ are shown in Figure~\ref{fig convex u set}. Each axis corresponds to one dimension of the uncertain parameter (i.e., $\tilde{y}_1, \tilde{y}_2$). The dashed lines correspond to the axis-aligned minimum bounding set $\cH$. The cyan shaded region represents an example of the admissible intervals for the uncertainty set which is not equivalent to the area formed by $[\buy,\boy]$, since $[\buy,\boy]$ may not be fully admissible and may lead to an ineffective budget of uncertainty. Finally, the orange shape inside the uncertainty set corresponds to the boundaries of the effective uncertainty set which is a convex set within the admissible box and its shape can be derived depending on the shape of the uncertainty set $\cY$. For example, in the polyhedral case with an initial box uncertainty, the effective uncertainty is also a box, as implied by Proposition \ref{prop worst-case realization interval}. However, the effective uncertainty for an ellipsoidal uncertainty set is a new ellipsoid as shown in Appendix \ref{Appendix ellipsoid}, Figure \ref{fig: eff ellipsoid}. 

{\color{black}
}
\begin{figure}[htbp]
    \centering
\includegraphics[width=0.96\textwidth]{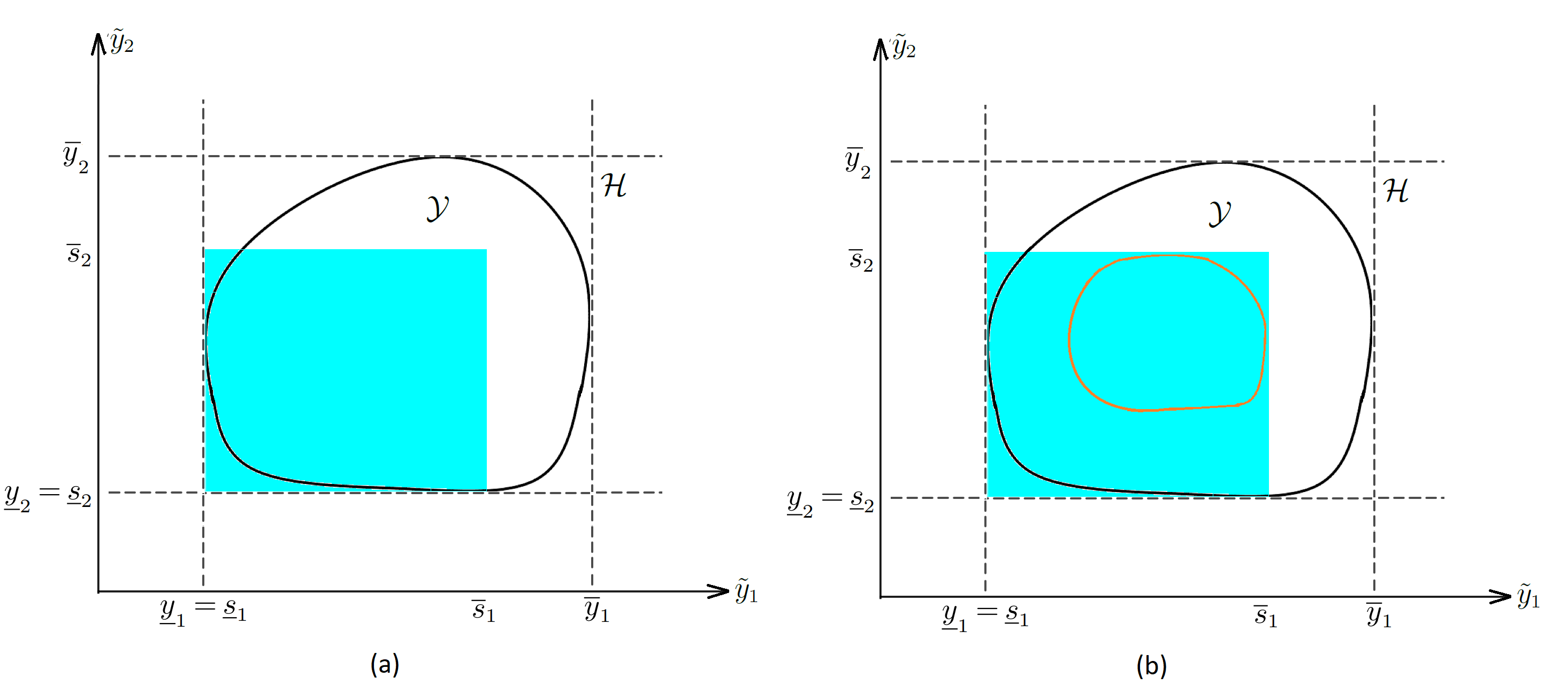}
\caption{(a) The uncertainty set $\cY$ and its admissible area (shaded area). (b) The effective uncertainty set of $\cY$ (shown in orange). Dotted lines in both figures correspond to hypercube $\cH$ containing $\cY$.  
    }
    \label{fig convex u set}
\end{figure}

Based on Remark~\ref{remark 1}, the admissible interval can be obtained using the optimization problem~\ref{sub} with no restrictions on the uncertain parameters but the lower and upper bounds (i.e., constraints \eqref{sub-d}-\eqref{sub-e}). In the general case, this problem corresponds to finding the admissible interval within $\cH$, the axis-aligned minimum bounding box of $\cY$. Thus, Remark \ref{remark 1} and Proposition \ref{prop - admissible interval} apply to a general case regardless of the shape of the underlying uncertainty set. Finally, we need to find the effective budgeted uncertainty set $\cU^{E}$ to only consider effective deviations that can potentially impact the solution. The specific steps of this general procedure is straightforward and requires the analysis the structure of the uncertainty set $\cU$ and its admissible part, $\cU^A$. As expected, the complexity of the problem depends on the shape of the underlying uncertainty set $\cY$ which affects the formulation of its robust counterpart. In Appendix~\ref{Appendix ellipsoid}, we provide details of these steps for the case of ellipsoidal uncertainty. 

} 

\section{Numerical Results} \label{results}
In this section, we examine the performance of the proposed robust model on an example of the previously-described multi-period SCED problem with wind uncertainty. The SCED problem can be presented in the form of model $\M$ where vectors 
$\bx$ and $\by$ correspond to the power generation of conventional generators and wind turbines, respectively. 
The objective function aims to minimize to the total operational cost, i.e., the generation cost of conventional generators plus a penalty cost associated with non-utilized wind power, under the worst-case scenario of wind availability over a given time horizon. Constraint \eqref{M2-b} captures the supply-demand balance equations and operational limits (e.g., limits of transmission lines, generators, and reserve requirements) of the power system. Finally, the upper bound on the available wind power is captured in constraint \eqref{M2-d}, where the utilized wind power is less than or equal to the available amount of wind power.

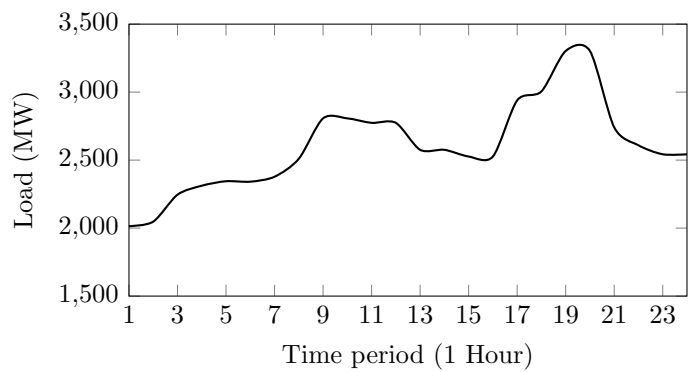
\begin{figure}[htbp] 
    \centering
  \begin{tikzpicture}[font=\small]
\pgfplotsset{width=9cm, height=5.2cm}
  \begin{axis}[
   xmin=1,
    xmax=24,
    xtick={1,3,...,23},
    xlabel= Time period (1 Hour),
    ymin=1500,
    ymax=3500,
    ytick={1500,2000,...,3500},
    ylabel={Load (MW)}
  ] 
   \addplot [no marks, black,thick, smooth,name path=load] plot coordinates {
(1,2014.464)   (2,2047.488) (3,2245.632)  (4,2311.68)    (5,2344.704)   (6,2341.4016)   (7,2377.728) (8,2509.824)   (9,2807.04)   (10,2807.04)  (11,2774.016) (12,2774.016)  (13,2575.872)   (14,2575.872)   (15,2526.336)   (16,2526.336)  (17,2939.136)   (18,3005.184)   (19,3302.4)   (20,3302.4)   (21,2740.992)   (22,2608.896)   (23,2542.848)  (24,2542.848)
   };
  \end{axis}
\end{tikzpicture}
    \caption{Load profile over a 24-hour time horizon}
    \label{load profile}
\end{figure}

\begin{figure}[t]
\hspace{-0.5cm}
  \centering
  \medskip
  \begin{subfigure}[t]{.5\linewidth}
    \centering
\begin{tikzpicture}[font=\small]
\pgfplotsset{width=8cm, height=5cm}
\begin{axis}[legend image post style={scale=0.5},legend style={{font=\fontsize{6.7}{6.7}\selectfont},
at={(0.0,0.00)}, anchor=south west},
   xmin=1,
    xmax=24,
    xtick={1,3,...,23},
    xlabel= Time period (1 Hour),
    ymin=100,
    ymax=640,
    ytick={100,250,...,700},
    ylabel={Power output (MW)}
  ] 
   \addplot [no marks, black,dashed,thick, smooth,name path=Upper limit of predicted wind power] plot coordinates {
(1,378.44)   (2,380.98) (3,386.54)  (4,377.1)    (5,384)   (6,421.24)   (7,444.84) (8,493.44)   (9,523.22)   (10,538.68)  (11,554.2) (12,571.98)  (13,564.4)   (14,558.36)   (15,548.42)   (16,527.78)  (17,514.62)   (18,480.26)   (19,478.24)   (20,490.78)   (21,490)   (22,485)   (23,490)  (24,500)
   };
%
   \addplot [no marks, red,dashed,thick, smooth,name path=Lower limit of predicted wind power] plot coordinates {
(1,321.56)   (2,309.02) (3,313.46)  (4,302.9)    (5,306)   (6,338.76)   (7,355.16) (8,396.56)   (9,416.78)   (10,421.32)  (11,425.8) (12,428.02)  (13,415.6)   (14,381.64)   (15,351.58)   (16,312.22)  (17,285.38)   (18,219.74)   (19,201.76)   (20,179.22)   (21,170) (22,165)   (23,130)  (24,100)};
   \addplot [no marks, brown ,thick, smooth,name path= Predicted wind power] plot coordinates {
(1,350)   (2,345) (3,350)  (4,340)    (5,345)   (6,380)   (7,400) (8,445)   (9,470)   (10,480)  (11,490) (12,500)  (13,490)   (14,470)   (15,450)   (16,420)  (17,400)   (18,350)   (19,340)   (20,335)   (21,330)   (22,325)   (23,310)  (24,300)};
\legend{Upper limit of wind power, Lower limit of wind power, Predicted wind power}
  \end{axis}
\end{tikzpicture}
\caption{Wind profile \#1}
  \end{subfigure}  
  \begin{subfigure}[t]{.5\linewidth}
    \centering
    \begin{tikzpicture}[font=\small]
\pgfplotsset{width=8cm, height=5cm}
\begin{axis}[
   xmin=1,
    xmax=24,
    xtick={1,3,...,23},
    xlabel= Time period (1 Hour),
    ymin=50,
    ymax=350,
    ytick={50,150,...,350},
    ylabel={Power output (MW)}
  ] 
   \addplot [no marks, black,dashed,thick, smooth,name path=Upper limit of predicted wind power] plot coordinates {
(1,96.34)   (2,107.6) (3,132.804) (4,142.984) (5,149.112) (6,144.436)   (7,160.312)(8,171.956)   (9,163.324)(10,180.544)  (11,193.208)(12,200.2)  (13,208.944)(14,230.924)   (15,263.732)(16,287.728)  (17,307.42)(18,315.588)   (19,309.296)(20,304.528)   (21,319.34)(22,324.524)   (23,338.676)  (24,345.588)
   };
%
   \addplot [no marks, red,dashed,thick, smooth,name path=Lower limit of predicted wind power] plot coordinates {
(1,73.66)   (2,82.4) (3,107.196)  (4,117.016)    (5,120.888)   (6,115.564)   (7,129.688) (8,138.044)   (9,126.676)   (10,139.456)  (11,146.792) (12,149.8)  (13,151.056)   (14,169.076)   (15,196.268)   (16,212.272)  (17,222.58)   (18,224.412)   (19,210.704)   (20,195.472)   (21,200.66) (22,195.476)   (23,201.324)  (24,194.412)};
   \addplot [no marks, brown ,thick, smooth,name path= Predicted wind power] plot coordinates {
(1,85)   (2,95) (3,120)  (4,130)    (5,135)   (6,130)   (7,145) (8,155)   (9,145)   (10,160)  (11,170) (12,175)  (13,180)   (14,200)   (15,230)   (16,250)  (17,265)   (18,270)   (19,260)   (20,250)   (21,260)   (22,260)   (23,270)  (24,270)};
  \end{axis}
\end{tikzpicture}
\caption{Wind profile \#2}
 \end{subfigure}
 \caption{Nominal value (predicted), upper limit, and lower limits of wind power}
 \label{Wind profile}
\end{figure}
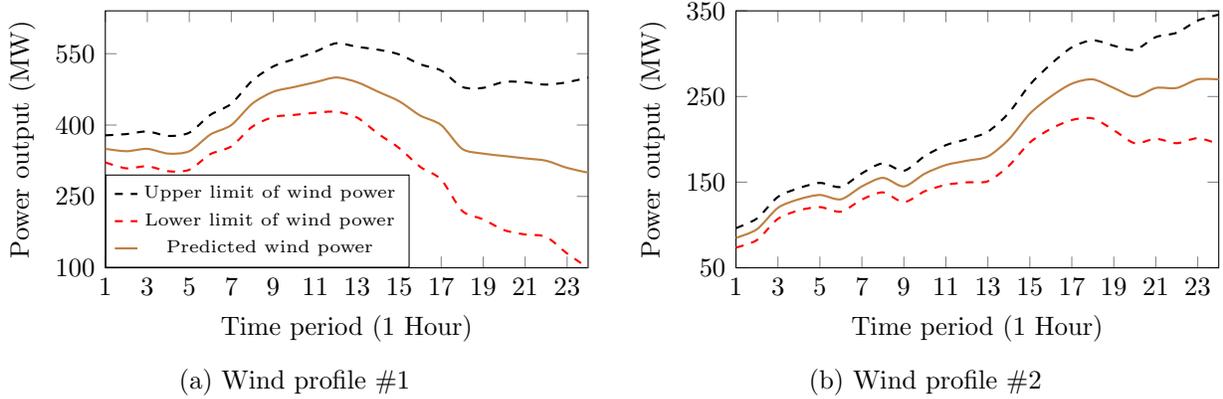

We use an IEEE reliability test system (RTS) with 24 buses and add multiple wind farms to perform numerical testing of our methodology. The test system is shown in Appendix~\ref{rts24}, and { The detailed data and parameters of the system, such as generation capacities, ramping rates, and line capacities, can be found in \cite{grigg1999ieee}. The generation cost data ($\bc_1$) are derived from the 
\cite{Washington}, and the penalty cost associated with wind power curtailment ($\bc_2$) is set to the maximum generation cost of conventional generators at 29.7 $\$$/MW \citep{li2016adaptive}, such that wind power units have dispatch priority over conventional generators.}
The hourly load (demand) is shown in Figure~\ref{load profile}. 
We consider a total of four wind farms in the system, two with nominal wind profile \#1 added to buses 13 and 16 $-$ denoted as wind farms A and B respectively $-$ and another two with nominal wind profile \#2 added to buses 18 and 8 $-$ denoted as wind farms C and D respectively. These nominal wind profiles and their corresponding uncertainty sets (upper and lower bounds) on the available wind power for each time period are shown in Figure~\ref{Wind profile}.  

We solve a day-ahead optimization model for the SCED problem with a 24-hour time horizon. The budget of uncertainty is defined for each one-hour time period across all wind farms. A number of power dispatch studies in the literature neglect the impacts of the unutilized wind power, i.e., the wind curtailment, in the optimal solution and assume that the nominal wind power can be fully utilized regardless of how volatile and large the available wind is \citep{wu2014robust}. We call this approach the ``naive approach'' where instead of optimizing over $\by$, we assume  $\by=\byhat$ is given. For comparison, we will use three variants of modeling the SCED optimization problem with wind power curtailment in this section: (i) the deterministic model {(formulation \eqref{M2})}, which does not account for wind uncertainty, and naively assumes the entire wind can be absorbed (ii) the conventional robust model {(formulation \eqref{M2 with Budget})}, which considers the absolute worst-case of wind uncertainty,
and (iii) the proposed two-stage robust model {(formulations \eqref{sub} and \eqref{master problem - stage II})} which considers the effective worst-case only. 
Detailed formulations of these three models for the SCED problem are provided in~\ref{appendix SCED}. 
For further details on the robust formulation of the SCED problem, the reader is referred to \cite{dehghani2019robust}. All models 
were solved to optimality using the IBM ILOG CPLEX solver with a solution time of a fraction of one second in all instances. 

%

In the rest of this section, we first perform day-ahead analyses where the uncertainty set on the available wind power for the next 24 hours is used to solve the SCED problem before knowing the actual realization of the uncertain parameter. We explicitly show all potential cases of the admissible interval, demonstrate the performance of the proposed approach in terms of the effectiveness of the budget of uncertainty, and compare the price of robustness in the conventional and proposed budget approaches.
Next, to verify the reliability and cost efficiency of the day-ahead solutions of the proposed approach, we compare the solutions of each approach with a  
“prescient” solution 
in which we hypothetically assume we have the perfect information of future wind power availability when planning in advance.
\subsection{Admissible Wind Power Interval} \label{results: admissible}
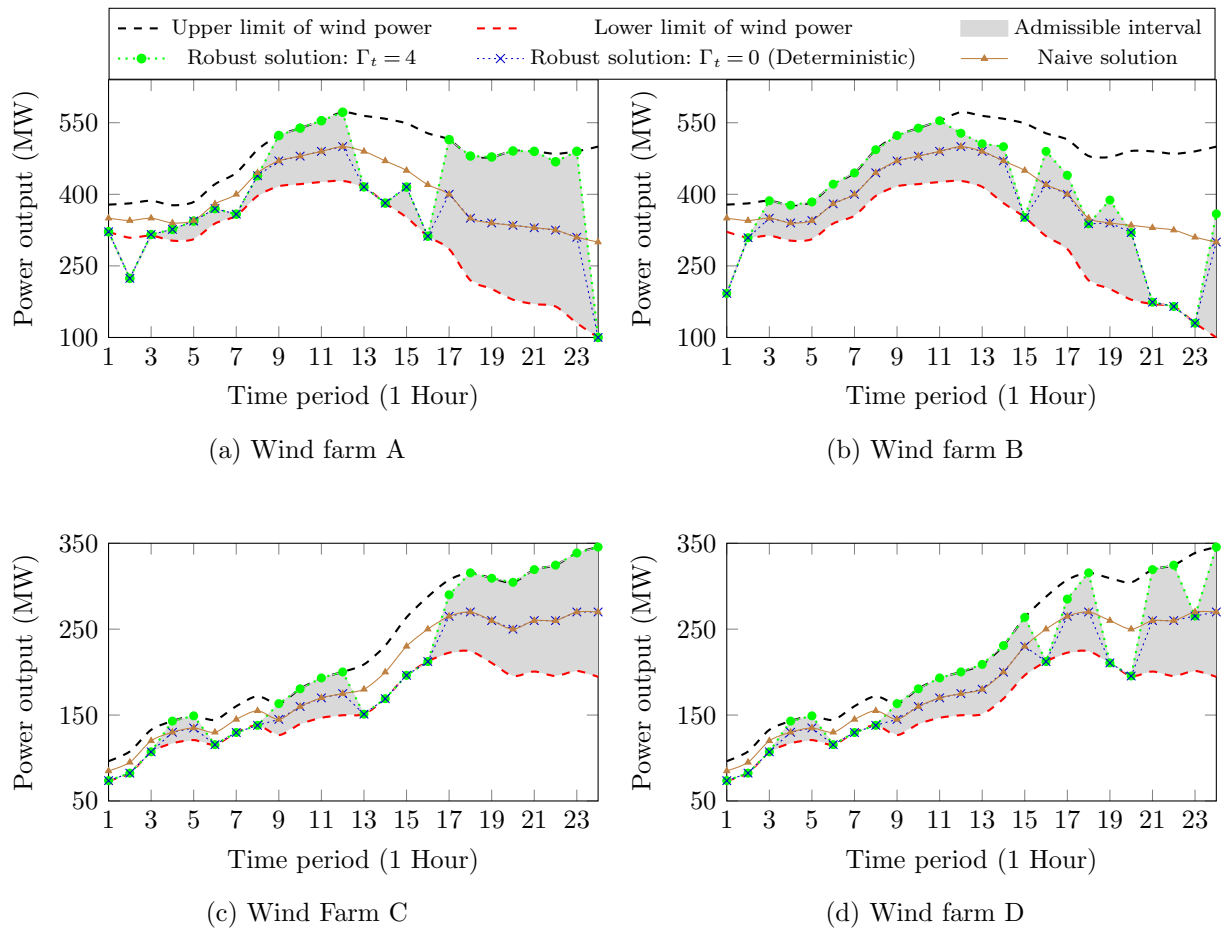
\begin{figure}[!t]
\centering
\medskip
\begin{subfigure}[t]{0.49\linewidth}
\begin{tikzpicture}[font=\small]
\pgfplotsset{width=1\textwidth, height=0.22\textheight}
\begin{axis}[legend columns=3, legend image post style={scale=1.04},
legend style={/tikz/every even column/.append style={column sep=0.51cm},{font=\fontsize{8}{8}\selectfont}, at={(0.0,1.0)}, anchor=south west},
   xmin=1,
    xmax=24,
    xtick={1,3,...,23},
    xlabel= Time period (1 Hour),
    ymin=100,
    ymax=640,
    ytick={100,250,...,700},
    ylabel={Power output (MW)}
  ] 
   \addplot [no marks, black,dashed,thick, smooth,name path=Upper limit of predicted wind power] plot coordinates {
(1,378.44)   (2,380.98) (3,386.54)  (4,377.1)    (5,384)   (6,421.24)   (7,444.84) (8,493.44)   (9,523.22)   (10,538.68)  (11,554.2) (12,571.98)  (13,564.4)   (14,558.36)   (15,548.42)   (16,527.78)  (17,514.62)   (18,480.26)   (19,478.24)   (20,490.78)   (21,490)   (22,485)   (23,490)  (24,500)
   };
%
   \addplot [no marks, red,dashed,thick, smooth,name path=Lower limit of predicted wind power] plot coordinates {
(1,321.56)   (2,309.02) (3,313.46)  (4,302.9)    (5,306)   (6,338.76)   (7,355.16) (8,396.56)   (9,416.78)   (10,421.32)  (11,425.8) (12,428.02)  (13,415.6)   (14,381.64)   (15,351.58)   (16,312.22)  (17,285.38)   (18,219.74)   (19,201.76)   (20,179.22)   (21,170) (22,165)   (23,130)  (24,100)};

\addplot [no marks,gray!30, thin , name path= zerohelp] plot coordinates {
(1,1)  (2,1)};
\addplot [no marks,gray!30, thin , name path= onehelp] plot coordinates {
(1,0)  (2,0)};
\addplot[gray!30] fill between[ of = zerohelp and onehelp];

   \addplot [no marks,gray!30, thin , name path= Upper limit of admissible wind power interval] plot coordinates {
(1,321.56)   (2,223.98) (3,315.748)  (4,326.572)    (5,343.896)   (6,370.104)   (7,358.384) (8,438.992)   (9,523.22)   (10,538.68)  (11,554.2) (12,571.98)  (13,415.6)   (14,381.64)   (15,415.18)   (16,312.22)  (17,514.62)   (18,480.26)   (19,478.24)   (20,490.78)   (21,490)   (22,468.237)   (23,490)  (24,100)};
%
   \addplot [no marks, dash pattern={on 0.0075pt off 100.4pt}, name path= Lower limit of admissible wind power interval] plot coordinates {
(1,321.56)   (2,223.98) (3,313.46)  (4,302.9)    (5,306)   (6,338.76)   (7,355.16) (8,396.56)   (9,416.78)   (10,421.32)  (11,425.8) (12,428.02)  (13,415.6)   (14,381.64)   (15,351.58)   (16,312.22)  (17,285.38)   (18,219.74)   (19,201.76)   (20,179.22)   (21,170)   (22,165)   (23,130)  (24,100)};

\addplot [no marks, gray!30, thick, name path= help Lower limit of admissible] plot coordinates {
(1,321.56)   (2,223.98) (3,313.46)  };
   \addplot [mark=*,mark options={scale=0.7,solid}, green!100, dotted, thick, name path= Robust solution full budget] plot coordinates {
(1,321.56)   (2,223.98) (3,315.748)  (4,326.572)    (5,343.896)   (6,370.104)   (7,358.384) (8,438.992)   (9,523.22)   (10,538.68)  (11,554.2) (12,571.98)  (13,415.6)   (14,381.64)   (15,415.18)   (16,312.22)  (17,514.62)   (18,480.26)   (19,478.24)   (20,490.78)   (21,490)   (22,468.237)   (23,490)  (24,100)};
%
   \addplot [mark=x, mark options={solid, scale=1.15}, blue!80!black, dash pattern={on 0.75pt off 1.4pt}, name path= Robust solution half budget] plot coordinates {
(1,321.56)   (2,223.98) (3,315.748)  (4,326.572)    (5,343.896)   (6,370.104)   (7,358.384) (8,438.992)   (9,470)   (10,480)  (11,490) (12,500)  (13,415.6)   (14,381.64)   (15,415.18)   (16,312.22)  (17,400)   (18,350)   (19,340)   (20,335)   (21,330)   (22,325)   (23,310)  (24,100)};

\addplot [mark=triangle*,mark options={scale=0.65,solid}, brown , smooth,name path= Predicted wind power] plot coordinates {
(1,350)   (2,345) (3,350)  (4,340)    (5,345)   (6,380)   (7,400) (8,445)   (9,470)   (10,480)  (11,490) (12,500)  (13,490)   (14,470)   (15,450)   (16,420)  (17,400)   (18,350)   (19,340)   (20,335)   (21,330)   (22,325)   (23,310)  (24,300)};
\addplot[gray!30] fill between[ 
    of = Upper limit of admissible wind power interval and Lower limit of admissible wind power interval];

\legend{Upper limit of wind power, Lower limit of wind power,,, Admissible interval,,,, Robust solution: $\Gamma_t=4$ ,Robust solution: $\Gamma_t=0$ (Deterministic) , Naive solution}
  \end{axis} 

\end{tikzpicture}
\caption{Wind farm A}
  \end{subfigure}
  \begin{subfigure}[t]{0.49\linewidth}
\begin{tikzpicture}[font=\small]
\pgfplotsset{width=1\textwidth, height=0.22\textheight}
\begin{axis}[legend columns=2, legend image post style={scale=0.95},
legend style={/tikz/every even column/.append style={column sep=0.8cm},{font=\fontsize{8}{8}\selectfont}, at={(0.0,1.0)}, anchor=south west},
   xmin=1,
    xmax=24,
    xtick={1,3,...,23},
    xlabel= Time period (1 Hour),
    ymin=100,
    ymax=640,
    ytick={100,250,...,700},
    ylabel={Power output (MW)}
  ] 
   \addplot [no marks, black,dashed,thick, smooth,name path=Upper limit of predicted wind power] plot coordinates {
(1,378.44)   (2,380.98) (3,386.54)  (4,377.1)    (5,384)   (6,421.24)   (7,444.84) (8,493.44)   (9,523.22)   (10,538.68)  (11,554.2) (12,571.98)  (13,564.4)   (14,558.36)   (15,548.42)   (16,527.78)  (17,514.62)   (18,480.26)   (19,478.24)   (20,490.78)   (21,490)   (22,485)   (23,490)  (24,500)
   };
%
   \addplot [no marks, red,dashed,thick, smooth,name path=Lower limit of predicted wind power] plot coordinates {
(1,321.56)   (2,309.02) (3,313.46)  (4,302.9)    (5,306)   (6,338.76)   (7,355.16) (8,396.56)   (9,416.78)   (10,421.32)  (11,425.8) (12,428.02)  (13,415.6)   (14,381.64)   (15,351.58)   (16,312.22)  (17,285.38)   (18,219.74)   (19,201.76)   (20,179.22)   (21,170) (22,165)   (23,130)  (24,100)};
%
   \addplot [no marks, gray!30, thin,name path= Upper limit of admissible wind power interval] plot coordinates {
(1,192.44)   (2,309.02) (3,386.54)  (4,377.1)    (5,384)   (6,421.24)   (7,444.84) (8,493.44)   (9,523.22)   (10,538.68)  (11,554.2) (12,527.98)  (13,505.88)   (14,499.84)   (15,351.58)   (16,489.996)  (17,440)   (18,338.201)   (19,388.021)   (20,319.937)   (21,173.973)   (22,165)   (23,130)  (24,358.981)};
%
   \addplot [no marks, gray!30, thick, name path= HELPLower limit of admissible wind power interval] plot coordinates {
(1,192.44)   (2,309.02)};

   \addplot [no marks, dash pattern={on 0.0075pt off 100.4pt}, name path= Lower limit of admissible wind power interval] plot coordinates {
(1,192.44)   (2,309.02) (3,313.46)  (4,302.9)    (5,306)   (6,338.76)   (7,355.16) (8,396.56)   (9,416.78)   (10,421.32)  (11,425.8) (12,428.02)  (13,415.6)   (14,381.64)   (15,351.58)   (16,312.22)  (17,285.38)   (18,219.74)   (19,201.76)   (20,179.22)   (21,170)   (22,165)   (23,130)  (24,100)};

   \addplot [mark=*,mark options={scale=0.65,solid}, green!100, dotted, thick, name path= Robust solution full budget] plot coordinates {
(1,192.44)   (2,309.02) (3,386.54)  (4,377.1)    (5,384)   (6,421.24)   (7,444.84) (8,493.44)   (9,523.22)   (10,538.68)  (11,554.2) (12,527.98)  (13,505.88)   (14,499.84)   (15,351.58)   (16,489.996)  (17,440)   (18,338.201)   (19,388.021)   (20,319.937)   (21,173.973)   (22,165)   (23,130)  (24,358.981)};
   \addplot [mark=x, mark options={solid, scale=1.15}, blue!80!black, dash pattern={on 0.75pt off 1.4pt}, name path= Robust solution half budget] plot coordinates {
(1,192.44)   (2,309.02) (3,350)  (4,340)    (5,345)   (6,380)   (7,400) (8,445)   (9,470)   (10,480)  (11,490) (12,500)  (13,490)   (14,470)   (15,351.58)   (16,420)  (17,400)   (18,338.201)   (19,340)   (20,319.937)   (21,173.973)   (22,165)   (23,130)  (24,300)};

\addplot [mark=triangle*,mark options={scale=0.65,solid}, brown , smooth,name path= Predicted wind power] plot coordinates {
(1,350)   (2,345) (3,350)  (4,340)    (5,345)   (6,380)   (7,400) (8,445)   (9,470)   (10,480)  (11,490) (12,500)  (13,490)   (14,470)   (15,450)   (16,420)  (17,400)   (18,350)   (19,340)   (20,335)   (21,330)   (22,325)   (23,310)  (24,300)};
\addplot[gray!30] fill between[ 
    of = Upper limit of admissible wind power interval and Lower limit of admissible wind power interval];
  \end{axis}
\end{tikzpicture}
\caption{Wind farm B}
  \end{subfigure}
  
  \vspace{2em}
  \begin{subfigure}[b]{0.49\linewidth}
\begin{tikzpicture}[font=\small]
\pgfplotsset{width=1\textwidth, height=0.22\textheight}
\begin{axis}[legend columns=2, legend image post style={scale=0.95},
legend style={/tikz/every even column/.append style={column sep=0.8cm},{font=\fontsize{8}{8}\selectfont}, at={(0.0,1.0)}, anchor=south west},
   xmin=1,
    xmax=24,
    xtick={1,3,...,23},
    xlabel= Time period (1 Hour),
    ymin=50,
    ymax=350,
    ytick={50,150,...,350},
    ylabel={Power output (MW)}
  ] 
   \addplot [no marks, black,dashed,thick, smooth,name path=Upper limit of predicted wind power] plot coordinates {
(1,96.34)   (2,107.6) (3,132.804) (4,142.984) (5,149.112) (6,144.436)   (7,160.312)(8,171.956)   (9,163.324)(10,180.544)  (11,193.208)(12,200.2)  (13,208.944)(14,230.924)   (15,263.732)(16,287.728)  (17,307.42)(18,315.588)   (19,309.296)(20,304.528)   (21,319.34)(22,324.524)   (23,338.676)  (24,345.588)
   };
%
   \addplot [no marks, red,dashed,thick, smooth,name path=Lower limit of predicted wind power] plot coordinates {
(1,73.66)   (2,82.4) (3,107.196)  (4,117.016)    (5,120.888)   (6,115.564)   (7,129.688) (8,138.044)   (9,126.676)   (10,139.456)  (11,146.792) (12,149.8)  (13,151.056)   (14,169.076)   (15,196.268)   (16,212.272)  (17,222.58)   (18,224.412)   (19,210.704)   (20,195.472)   (21,200.66) (22,195.476)   (23,201.324)  (24,194.412)};
%
   \addplot [no marks, gray!30, thin, name path= Upper limit of admissible wind power interval] plot coordinates {
(1,73.66)   (2,82.4) (3,107.196)  (4,142.984)    (5,149.112)   (6,115.564)   (7,129.688) (8,138.044)   (9,163.324)   (10,180.544)  (11,193.208) (12,200.2)  (13,151.056)   (14,169.076)   (15,196.268)   (16,212.272)  (17,290)   (18,315.588)   (19,309.296)   (20,304.528)  (21,319.34)   (22,324.524)   (23,338.676) (24,345.588)};
%
   \addplot [no marks, dash pattern={on 0.0075pt off 100.4pt}, name path= Lower limit of admissible wind power interval] plot coordinates {
(1,73.66)   (2,82.4) (3,107.196)  (4,117.016)    (5,120.888)   (6,115.564)   (7,129.688) (8,138.044)   (9,126.676)   (10,139.456)  (11,146.792) (12,149.8)  (13,151.056)   (14,169.076)   (15,196.268)   (16,212.272)  (17,222.58)   (18,224.412)   (19,210.704)   (20,195.472)   (21,200.66) (22,195.476)   (23,201.324)  (24,194.412)};

   \addplot [mark=*,mark options={scale=0.65,solid}, green!100, dotted, thick, name path= Robust solution full budget] plot coordinates {
(1,73.66)   (2,82.4) (3,107.196)  (4,142.984)    (5,149.112)   (6,115.564)   (7,129.688) (8,138.044)   (9,163.324)   (10,180.544)  (11,193.208) (12,200.2)  (13,151.056)   (14,169.076)   (15,196.268)   (16,212.272)  (17,290)   (18,315.588)   (19,309.296)   (20,304.528)  (21,319.34)   (22,324.524)   (23,338.676) (24,345.588)};
   \addplot [mark=x, mark options={solid, scale=1.15}, blue!80!black, dash pattern={on 0.75pt off 1.4pt}, name path= Robust solution half budget] plot coordinates {
(1,73.66)   (2,82.4) (3,107.196)  (4,130)    (5,135)   (6,115.564)   (7,129.688) (8,138.044)   (9,145)   (10,160)  (11,170) (12,175)  (13,151.056)   (14,169.076)   (15,196.268)   (16,212.272)  (17,265)   (18,270)   (19,260)   (20,250)   (21,260)   (22,260)   (23,270)  (24,270)};

\addplot [mark=triangle*,mark options={scale=0.65}, brown , smooth,name path= Predicted wind power] plot coordinates {
(1,85)   (2,95) (3,120)  (4,130)    (5,135)   (6,130)   (7,145) (8,155)   (9,145)   (10,160)  (11,170) (12,175)  (13,180)   (14,200)   (15,230)   (16,250)  (17,265)   (18,270)   (19,260)   (20,250)   (21,260)   (22,260)   (23,270)  (24,270)};
\addplot[gray!30] fill between[ 
    of = Upper limit of admissible wind power interval and Lower limit of admissible wind power interval];
  \end{axis}
\end{tikzpicture}
\caption{Wind Farm C}
  \end{subfigure}
    \begin{subfigure}[b]{0.49\linewidth}
\begin{tikzpicture}[font=\small]
\pgfplotsset{width=1\textwidth, height=0.22\textheight}
\begin{axis}[legend columns=2, legend image post style={scale=0.95},
legend style={/tikz/every even column/.append style={column sep=0.8cm},{font=\fontsize{8}{8}\selectfont}, at={(0.0,1.0)}, anchor=south west},
   xmin=1,
    xmax=24,
    xtick={1,3,...,23},
    xlabel= Time period (1 Hour),
    ymin=50,
    ymax=350,
    ytick={50,150,...,350},
    ylabel={Power output (MW)}
  ] 
   \addplot [no marks, black,dashed,thick, smooth,name path=Upper limit of predicted wind power] plot coordinates {
(1,96.34)   (2,107.6) (3,132.804) (4,142.984) (5,149.112) (6,144.436)   (7,160.312)(8,171.956)   (9,163.324)(10,180.544)  (11,193.208)(12,200.2)  (13,208.944)(14,230.924)   (15,263.732)(16,287.728)  (17,307.42)(18,315.588)   (19,309.296)(20,304.528)   (21,319.34)(22,324.524)   (23,338.676)  (24,345.588)
   };
%
   \addplot [no marks, red,dashed,thick, smooth,name path=Lower limit of predicted wind power] plot coordinates {
(1,73.66)   (2,82.4) (3,107.196)  (4,117.016)    (5,120.888)   (6,115.564)   (7,129.688) (8,138.044)   (9,126.676)   (10,139.456)  (11,146.792) (12,149.8)  (13,151.056)   (14,169.076)   (15,196.268)   (16,212.272)  (17,222.58)   (18,224.412)   (19,210.704)   (20,195.472)   (21,200.66) (22,195.476)   (23,201.324)  (24,194.412)};
%
   \addplot [no marks, gray!30 , thin,name path= Upper limit of admissible wind power interval] plot coordinates {
(1,73.66)   (2,82.4) (3,107.196)  (4,142.984)    (5,149.112)   (6,115.564)   (7,129.688) (8,138.044)   (9,163.324)   (10,180.544)  (11,193.208) (12,200.2)  (13,208.944)   (14,230.924)   (15,263.732)   (16,212.272)  (17,285)   (18,315.588)   (19,210.704)   (20,195.472)   (21,319.34) (22,324.524)   (23,265.305)  (24,345.588)};
%
   \addplot [no marks,dash pattern={on 0.0075pt off 100.4pt}, name path= Lower limit of admissible wind power interval] plot coordinates {
(1,73.66)   (2,82.4) (3,107.196)  (4,117.016)    (5,120.888)   (6,115.564)   (7,129.688) (8,138.044)   (9,126.676)   (10,139.456)  (11,146.792) (12,149.8)  (13,151.056)   (14,169.076)   (15,196.268)   (16,212.272)  (17,222.58)   (18,224.412)   (19,210.704)   (20,195.472)   (21,200.66) (22,195.476)   (23,201.324)  (24,194.412)};

   \addplot [mark=*,mark options={scale=0.65,solid}, green!100, dotted, thick, name path= Robust solution full budget] plot coordinates {
(1,73.66)   (2,82.4) (3,107.196)  (4,142.984)    (5,149.112)   (6,115.564)   (7,129.688) (8,138.044)   (9,163.324)   (10,180.544)  (11,193.208) (12,200.2)  (13,208.944)   (14,230.924)   (15,263.732)   (16,212.272)  (17,285)   (18,315.588)   (19,210.704)   (20,195.472)   (21,319.34) (22,324.524)   (23,265.305)  (24,345.588)};
   \addplot [mark=x, mark options={solid, scale=1.15}, blue!80!black, dash pattern={on 0.75pt off 1.4pt}, name path= Robust solution half budget] plot coordinates {
(1,73.66)   (2,82.4) (3,107.196)  (4,130)    (5,135)   (6,115.564)   (7,129.688) (8,138.044)   (9,145)   (10,160)  (11,170) (12,175)  (13,180)   (14,200)   (15,230)   (16,212.272)  (17,265)   (18,270)   (19,210.704)   (20,195.472)   (21,260)   (22,260)   (23,265.305)  (24,270)};

\addplot [mark=triangle*,mark options={scale=0.65}, brown , smooth,name path= Predicted wind power] plot coordinates {
(1,85)   (2,95) (3,120)  (4,130)    (5,135)   (6,130)   (7,145) (8,155)   (9,145)   (10,160)  (11,170) (12,175)  (13,180)   (14,200)   (15,230)   (16,250)  (17,265)   (18,270)   (19,260)   (20,250)   (21,260)   (22,260)   (23,270)  (24,270)};
\addplot[gray!30] fill between[ 
    of = Upper limit of admissible wind power interval and Lower limit of admissible wind power interval];
  \end{axis}
\end{tikzpicture}
\caption{Wind farm D}
  \end{subfigure}
 \caption{Wind power output comparison between the proposed robust, deterministic, and naive approaches. Shaded regions correspond to the wind power admissible intervals.
 }
 \label{wind generation}
\end{figure}

{ Figure~\ref{wind generation} shows the admissible intervals for wind power for all periods across all wind farms, as well as the day-ahead solution of three different approaches. The two dashed lines in the Figure show the upper and lower bounds of the wind power uncertainty set, as shown in Figure~\ref{Wind profile}, and the middle line represents the nominal wind power and the naive solution that assumes all the available wind can be absorbed. The shaded regions show the admissible intervals for each time period, obtained by solving the first-stage problem in formulation \eqref{sub}. 
We can observe that different parts of the shaded regions for different periods correspond to examples of all four possible cases of the admissible interval mentioned in Proposition~\ref{prop - admissible interval} during this 24-hour time horizon. For instance, During periods 9 to 12, wind farm A has the admissible interval of case (a), where the available wind power can be entirely absorbed.} Cases (b) and (c) are observed in wind farm B during period 16 and wind farm A during period 15, respectively, in which a part of the available wind power cannot be utilized. Case (d) is observed in wind farm A during period 2, where the admissible wind power interval is outside of the uncertainty set.

Figure~\ref{wind generation} also shows the day-ahead robust solutions of the proposed approach using dotted lines with $\Gamma_t=4$ (full budget) and $\Gamma_t=0$ (no budget) in comparison with that of the naive approach. Note that the proposed robust solution with $\Gamma_t=0$ corresponds to the deterministic solution with no wind power uncertainty. From Figure~\ref{wind generation}, we observe that the proposed robust solutions lie within the admissible wind power interval (shaded region) and thus guarantee feasibility (system security) under any scenario of the actual wind power in real-time. {However, this is not the case for the naive approach and the corresponding solution may not be feasible under all wind power scenarios in real-time}. For example, during some periods (e.g., periods 1 and 2), even though the naive solution corresponds to more wind power utilization, the solution does not lie within the admissible interval, and thus it may violate the operational limits of the systems in real-time after knowing the actual realization of wind power since the system capabilities might not be sufficient to adjust the naive day-ahead solution. Furthermore, it can be observed that the effective worst-case scenario respects reliability and is always a part of the admissible interval (the upper bound of the admissible interval for $\Gamma=4$ { shown with green dotted line}) because it is obtained after taking the worst-case scenario into account. 
\subsection{Performance of The Proposed Approach versus The Conventional Budget Approach}\label{results: budget effectiveness}
\begin{figure}[t] 
    \centering
  \begin{tikzpicture}[font=\small]
\pgfplotsset{width=10cm, height=5.5cm}
  \begin{axis}[legend image post style={scale=0.6},legend style={{font=\fontsize{8}{8}\selectfont}, at={(0.0,1)}, anchor=north west},
tick label style={font=\fontsize{8.8}{8.8}\selectfont},
   xmin=0,
    xmax=4,
    xtick={0,0.5,...,4},
    xlabel= $\Gamma_{17}$,
    ymin=1330,
    ymax=1560,
    ytick={1330,1380,...,1530},
    ylabel={Utilized wind power (MW)}
  ] 
   \addplot [no marks, black,thick,name path=proposed approach] plot coordinates {
(0,1330) (0.5,1387.3) (1,1444.62) (1.35,1484) (1.5,1491.03) (2,1512.24) (2.47,1523) (2.5,1524)(3,1532) (3.5,1540)(3.59,1541)(4,1541)};
\addplot [no marks, red,dashed,thick, name path=conventional approach] plot coordinates {
(0,1330) (0.5,1387.3) (1,1444.62) (1.35,1484) (1.5,1484.737) (2,1484.737) (2.47,1495.6744) (2.5,1495.6744)(3,1495.6744) (3.5,1535)(3.59,1540)(4,1541)};
%
\addplot [no marks, ultra thin, gray!50, name path=h1] plot coordinates {(1.35,1484)(1.35,0)};
\addplot [no marks,  ultra thin, gray!50, name path=h2] plot coordinates {(2,1484)(2,0)};
\addplot [no marks,  ultra thin, gray!50, name path=h3] plot coordinates {(2.46,1495)(2.46,0)};
\addplot [no marks,  ultra thin, gray!50, name path=h4] plot coordinates {(3,1495)(3,0)};

\legend{Proposed approach, Conventional approach}
\end{axis}
\end{tikzpicture}
\caption{Comparison of the effectiveness between the proposed approach and the conventional approach in terms of wind power utilization during time period 17.}
\label{utilized wind}
\end{figure}
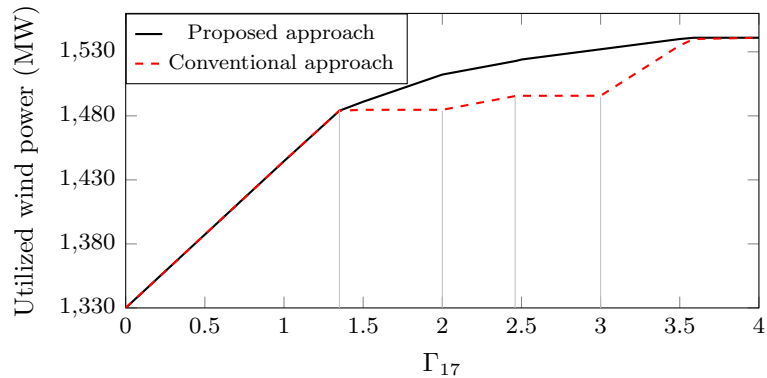

In this section, we compare the effectiveness of the proposed robust approach 
with that of the conventional approach 
for different levels of the budget of uncertainty. 
Recall that the budget of uncertainty is defined per time period. Here, we focus on comparing results for a sample time period (time period $t=17$), but we note that similar observations can be made for other time periods as well. Figure~\ref{utilized wind} shows the wind power utilization versus the budget of uncertainty for the sample time period. 
Overall, Figure~\ref{utilized wind} 
shows that the proposed approach results in larger wind power utilization compared to the conventional approach and adjusts the solution for various levels of the budget of uncertainty, as intuitively expected. In what follows, we explain the detailed reasoning behind this behaviour for $t=17$.

For values of $ \Gamma_{17} \in [0,1]$, we observe that both the conventional and proposed approaches result in the same wind power utilization, and an increase in the budget of uncertainty leads to higher wind power utilization. Here, the conventional approach uses the first unit of the budget for 
deviations of the uncertain parameter from the nominal value in wind farm A, whose initial uncertainty set corresponds to case (a) during time 17 and can be entirely utilized. Thus, the uncertain parameters are within the admissible interval for both approaches, and larger budgets of uncertainty correspond to  more wind power utilization without causing infeasibility. 

For $\Gamma_{17} \in (1,1.35]$, the uncertain parameter of the conventional approach corresponds to the admissible interval of wind farm B and results in the same solution as the proposed approach. On the other hand, for $\Gamma_{17} \in (1.35,2]$, the conventional approach leads to an ineffective budget that cannot be further utilized in wind farm B since it considers the corresponding absolute worst-case scenario as opposed to the admissible worst-case. Thus, the uncertain parameters corresponding to the available wind power lie outside of the admissible interval and cannot be utilized. Hence, changing the budget does not impact the robust solution anymore. 
 { On the contrary, the proposed robust approach leads to a less conservative solution with higher wind power utilization by removing the unnecessary protection for the ineffective parts of the uncertainty set. 
By only considering the effective worst-case scenario, it allocates the budget of uncertainty to a different wind farm for which the uncertain parameter is admissible.} We note that this corresponds to our definition of the effective worst-case and is not equivalent with the absolute worst-case in terms of cost. 

For $\Gamma_{17} \in (2,2.47]$ an increase in the budget of uncertainty results in higher wind power utilization in both approaches since the uncertain parameter corresponding to the power output of both wind farms are within the admissible interval of wind farm C. 
Similar arguments can be made for $\Gamma_{17} \in (2.47,3]$, $\Gamma_{17} \in (3, 3.59]$, and  $\Gamma_{17} \in (3.59, 4]$.  
Finally, for  $\Gamma_{17} =4$ (full budget), both approaches lead to 
the same overly-conservative solution, where the uncertain parameters have maximum deviations from their nominal values in all wind farms.

\subsection{Total Cost versus Budget of Uncertainty} \label{results: price of robustness}
\begin{figure}[t] 
    \centering
  \begin{tikzpicture}[font=\fontsize{9.5}{9.5}\selectfont]
\pgfplotsset{width=10cm, height=5.1cm}
  \begin{axis}[legend image post style={scale=0.6},legend style={{font=\fontsize{8}{8}\selectfont}, at={(0.0,1)}, anchor=north west},
   xmin=0,
    xmax=4,
    xtick={0,0.5,1,...,4},
    xlabel= $\Gamma_{t}$,
    ymin=810,
    ymax=885,
    ytick={810,830,...,900},
    ylabel={Total cost ($\times \$10^3$)}
  ] 
   \addplot [mark=triangle*, black,thick,name path=proposed approach] plot coordinates {
(0,818) (0.5,826) (1,838) (1.5,840) (2,842) (2.5,847) (3,856) (3.5,862)(4,868)};
\addplot [no marks, red,dashed,thick, name path=conventional approach] plot coordinates {

(0,818) (0.5,828) (1,845) (1.5,855) (2,858) (2.5,860) (3,862) (3.5,864)(4,868)};
%
\addplot [no marks, black,dashed,thick, name path=deterministic approach] plot coordinates {
(0,818) (4,818)};
\addplot[gray!30] fill between[ 
    of = conventional approach and proposed approach];
\legend{Proposed approach, Conventional approach, Deterministic approach}
  \end{axis}
\end{tikzpicture}
\caption[Short caption]{Total cost (objective function value) under various budgets of uncertainty.}
\label{price final}
\end{figure}
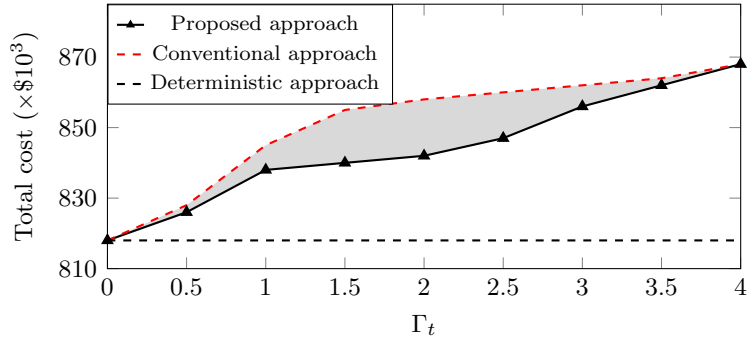

Figure~\ref{price final} shows the trade-off between the total operational cost and budget of uncertainty over the 24-hour time horizon of the day-ahead robust and deterministic solutions versus the budget of uncertainty. We note that these costs are the expected day-ahead cost and are not necessarily equal to the actual 
cost when the actual wind power realized. Later in this section, we will elaborate on the comparison of how close of an estimate of the actual cost each of these methods would provide
(Section~\ref{results: simulation}).

For the deterministic approach, the optimal operational cost 
is constant and does not depend on the budget of uncertainty. For the robust approaches, on the other hand, the operational cost consistently increases for higher budgets of uncertainty since a higher budget leads to a more conservative solution with higher day-ahead wind power curtailment during critical periods and this increases the objective function value due to high curtailment costs \citep{li2016adaptive}. This trade-off in the conventional budget approach is referred to as ``the price of robustness'' \citep{bertsimas2004price}, and it is based on the absolute worst-case scenario. We can also observe such a trade-off in the proposed approach. However, this trade-off accounts for the effective worst-case scenario only and corresponds to ``the effective price of robustness'' where the solution conservatism 
can be adjusted, and the operational cost is lower for $ 0 < \Gamma_t < 4$. { The shaded region demonstrates the difference between 
the proposed and conventional robust approaches in terms of their total cost. This reduction is cost is due to removing the unnecessary protection for the ineffective parts of the uncertainty set. By doing so, the budget of uncertainty can be allocated to the effective parts of the uncertainty set only, and the resulting robust solution is, therefore, less conservative.}

On the other hand, for $\Gamma_t=0$ (no budget) and $\Gamma_t=4$ (full budget), both robust approaches are equivalent since they both correspond to the deterministic and the traditional robust approach with over-conservative solutions, respectively.
It is to be mentioned that the actual realization of the wind power may still occur outside of the admissible region but when such a scenario is realized, it cannot be utilized, and considering the absolute worst-case scenario in the conventional budget approach leads to a conservative solution due to having a higher operational cost, as opposed to the effective worst-case scenario.

\subsection{Simulation Results and Cost Verification} \label{results: simulation}
One advantage of using a robust optimization methodology is that the resulting optimal solutions are often less sensitive to changes in the uncertain parameter compared to that of a deterministic approach and provides a better estimate of the realized cost, regardless of the realized scenario of the uncertain parameter. To demonstrate this advantage in our specific application, we 
compare the day-ahead solution ($\by$) of each approach 
with a \emph{prescient} solution which 
assumes that perfect information of $\by^{Actual}$ is known in advance and is obtained by solving the deterministic model with $\bty = \by^{Actual}$. 
{ { To obtain a set of simulated
wind scenarios, we generated randomized values of the actual wind power within interval $[\buy,\boy]$ for each time period and
calculated the scaled deviation \[\bz=\frac{|\by^{Actual}-\byhat|}{\boy-\byhat} \] from the nominal value for each period. This process was repeated until 100 scenarios are generated for each time interval such that $\sum_{i=1}^m z^+_i \le \Gamma_t$, for each of the different values of $\Gamma_t =$1, 2, 3, and 4. 
Note that for each time period, a scenario is a collection of randomized values for $\by^{Actual}$ of different wind farms such that the total deviations from the nominal values fits the budget constraint.}}
For each scenario, we solve the deterministic SCED problem 
with $\bty = \by^{Actual}$, and calculate the absolute difference between the total cost of the day-ahead and prescient solutions, denoted as $\Delta C$, for each of the three approaches. { This $\Delta C$ captures the wind power over-generation or under-generation costs with respect to $\by^{Actual}$. In our post-optimization real-time analysis, we consider a penalty cost of 50~\$/MWh for each unit of wind power over-generation or under-generation.} 

\begin{figure}[t]
\centering
\medskip
\begin{subfigure}[t]{0.31\linewidth}
\begin{tikzpicture}[font=\small]
\pgfplotsset{width=0.89\textwidth, height=0.28\textheight}
\begin{axis}[legend image post style={scale=0.8},legend style={{font=\fontsize{8}{8}\selectfont}, at={(0.02,0.98)}, anchor=north west},
   xmin=0.6,
    xmax=4.4,
    xtick={1,2,3,4},
    xlabel= $\Gamma_{t}$,
    ymin=5,
    ymax=45,
    ytick={5,15,...,45},
    ylabel={$\Delta C \,\,\, (\times \$10^3)$}] 
\addplot[scatter/classes={
a={mark=*,black}},
scatter,only marks,
scatter src=explicit symbolic] coordinates {(1,18.07)[a] (2,22)[a] (3,25.33)[a] (4,36)[a]};
\addplot [mark=-, black] plot coordinates {  (1,16) (1,19.42)
};
\addplot [mark=-, black] plot coordinates {  (2,24) (2,19)
};
\addplot [mark=-, black] plot coordinates {  (3,28.26) (3,21.50)
};
\addplot [mark=-, black] plot coordinates {  (4,43) (4,29)
};
\legend{Average}
\end{axis}
\end{tikzpicture}
\caption{Deterministic approach}
  \end{subfigure}
  \begin{subfigure}[t]{0.31\linewidth}
\begin{tikzpicture}[font=\small]
\pgfplotsset{width=0.89\textwidth, height=0.28\textheight}
\begin{axis}[legend image post style={scale=0.8},legend style={{font=\fontsize{8}{8}\selectfont}, at={(0.01,0.98)}, anchor=north west},
   xmin=0.6,
    xmax=4.4,
    xtick={1,2,3,4},
    xlabel= $\Gamma_{t}$,
   ymin=5,
    ymax=45,
    ytick={5,15,...,45},    
    ylabel={$\Delta C \,\,\, (\times \$10^3)$}] 
%
\addplot[scatter/classes={
a={mark=*,black}},
scatter,only marks,
scatter src=explicit symbolic] coordinates {(1,13.90)[a] (2,14.43)[a] (3,15.43)[a] (4,16.15)[a]};
\addplot [mark=-, black] plot coordinates {  (1,12) (1,15.3)
};
\addplot [mark=-, black] plot coordinates {  (2,17) (2,12.18)
};
\addplot [mark=-, black] plot coordinates {  (3,17.9) (3,12.2)
};
\addplot [mark=-, black] plot coordinates {  (4,12.71) (4,19.84)
};
\end{axis}
\end{tikzpicture}
\caption{Conventional approach}
  \end{subfigure}
  \begin{subfigure}[t]{0.31\linewidth}
\begin{tikzpicture}[font=\small]
\pgfplotsset{width=0.89\textwidth, height=0.28\textheight}
\begin{axis}[legend image post style={scale=0.8},legend style={{font=\fontsize{8}{8}\selectfont}, at={(0.01,0.98)}, anchor=north west},
   xmin=0.6,
    xmax=4.4,
    xtick={1,2,3,4},
    xlabel= $\Gamma_{t}$,
    ymin=5,
    ymax=45,
    ytick={5,15,...,45},
    ylabel={$\Delta C \,\,\, (\times \$10^3)$}] 
\addplot[scatter/classes={
a={mark=*,black}},
scatter,only marks,
scatter src=explicit symbolic] coordinates {(1,10.64)[a] (2,11.18)[a] (3,12.96)[a] (4,16.15)[a]};
\addplot [mark=-, black] plot coordinates {  (1,12.13) (1,9.2)
};
\addplot [mark=-, black] plot coordinates {  (2,9.23) (2,13.37)
};
\addplot [mark=-, black] plot coordinates {  (3,15.57) (3,10.24)
};
\addplot [mark=-, black] plot coordinates {  (4,19.84) (4,12.71)
};
\end{axis}
\end{tikzpicture}
\caption{Proposed approach}
  \end{subfigure}
\caption{The sensitivity of day-ahead solutions to the prescient solution with perfect realizations of wind. }
\label{simulationfinal}
\end{figure}
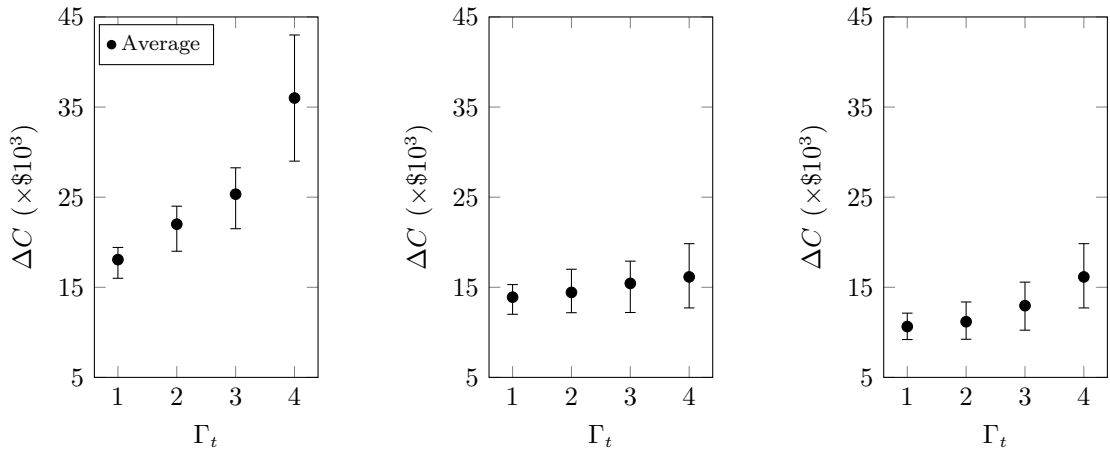

Figure~\ref{simulationfinal} compares the performance of the three approaches: 
(i) the deterministic approach, (ii) the conventional budget approach, and (iii) the proposed budget approach, in terms of $\Delta C$ which shows the sensitivity of day-ahead solutions to simulated realizations of wind power and provides a basis to compare which model gives a better estimate of the prescient solution.
We can observe that the deterministic approach corresponds to the largest $\Delta C$ values since this approach does not take into account the wind power uncertainty, and therefore the day-ahead scheduled (expected) cost is much lower than the prescient (actual) cost. 
On the other hand, the conventional and proposed robust approaches corresponds to smaller $\Delta C$ values as they account for uncertainty, and hence their day-ahead solutions are closer to the prescient solutions. 

The proposed robust approach has the smallest $\Delta C$ since it cuts off the ineffective part of the original uncertainty set and therefore the effective worst-case scenario falls within the admissible interval (similar to the prescient solution) 
to meet the operational limits of the system.
Thus, the proposed approach corresponds to a robust solution closer to the prescient solution. In other words, since the effective worst-case scenario corresponds to a smaller interval, i.e., $[\shj, \soj]$ as opposed to $[\yhj, \yoj]$, it is less sensitive to changes in simulations.
Particularly, for $\Gamma_t=1$, the average $\Delta C$ for the proposed robust approach is approximately 25\% and 40\% less than those of the conventional robust and deterministic approaches, respectively. This is a desirable outcome since it helps planners to have a much better estimate of the actual costs when performing day-ahead planning.  

\section{Conclusions} \label{conclusion}
{ 
In this paper, we proposed a robust optimization approach for considering an effective budget of uncertainty in a class of problems with 
right-hand side (RHS) uncertainty. 
We discussed limitations of the conventional budget-of uncertainty approach in such problems where adjusting the budget may not have an explicit effect on the solution. We 
proposed a two-stage robust approach that identifies the ineffective parts of the uncertainty set,  removes any unnecessary protection, 
and hence, produces less conservative solutions. We provide a one-to-one mapping between the budget parameter in our approach and that of the original budget-of-uncertainty approach. 
The proposed budget parameter approach allows the managers to make more informed 
decisions on the level of conservatism 
of the robust solution. 

We demonstrated the practical merits of the proposed approach in a security-constrained economic dispatch (SCED) problem with wind power uncertainty and numerically compared our results with those of the deterministic approach and the conventional robust approach. We also discussed other potential applications along with small numerical examples to demonstrate the applicability in each case.
%
An area of further research is to extend the proposed approach for other real-world applications with RHS uncertainty $-$ such as project management, scheduling, and telecommunication problems $-$ to more clearly understand the trade-off between the level of conservatism and budget of uncertainty for such problems. 
}

\bibliography{sample}
\bibliographystyle{apa} 

\begin{APPENDICES}
\section{Proof of Proposition \ref{prop first} }\label{Appendix proof}
In $\MBm$, constraint \eqref{M2-d budget} corresponds to $\by \le \boy$ for a full budget, which is satisfied for any value $\by$ within the admissible interval. Also,
variable $\by \in [\bus, \, \bos]
$ 
can be written as $\by = \bus+\br \odot (\bos-\bus)$, where $\bzero \le \br \le \bone$. Thus, for each row $i$, constraint \eqref{M2-b budget} can be reformulated as
$\bA_i\bx+\bB_i\by \le \bA_i\bx+\bB_i\bus+\bB_i(\br\odot(\bos-\bus))\le g_i$. 
To find the maximum feasible $\by$, the following constraint must be satisfied for each row $i$:
\begin{equation}\label{e5}
\bA_i\bx+ \bB_i\by \le \bA_i\bx + \bB_i\bus + \bm{\beta}_i (\br) \le g_i, \qquad \qquad \forall i \in \{1, \dots, n \} 
\end{equation}
where 
\begin{subequations}\label{e5-2}
 \begin{align} 
 \qquad \bm{\beta}_i (\br)= \max_{\br} & \quad 
 \bB_i \Big( \br\odot(\bos-\bus) \Big), \qquad \quad \forall i \in \{1, \dots, n \}  \\
 \text{s.t.} & \quad 
 \bzero \le \br \le \bone. & \label{e5-2b} 
    \end{align}
\end{subequations}
Considering $\balpha_i'$ as a column of dual vectors corresponding to constraint \eqref{e5-2b} and replacing the dual of problem \eqref{e5-2} into \eqref{e5}, we can recover formulation \eqref{M3 upper bound}.  

\vspace{1em}

{
\section{The case of Ellipsoidal Uncertainty} \label{Appendix ellipsoid}
An ellipsoidal uncertainty set considers the correlation between uncertain parameters and is 
relevant in certain applications \citep{kurtz2018robust}. 
Let $\bSigma \in \mathbb{R}^{n\times n}$ be a given positive definite matrix indicating the correlations. 
%
%
Let $\lambda_1, \dots, \lambda_m >0$ be the eigenvalues of matrix $\bSigma$, and $v_i, \dots, v_j$ be the corresponding unit eigenvectors that show the direction of the $j^\text{th}$ axis of the ellipsoid. The underlying uncertainty set without budget can be written as
\begin{equation}
\label{ellipsoide U set}
\cY_e = \Big\{ \hspace{0.051cm} \bty \in \mathbb{R}^{{m}} : \hspace{0.2cm}  ( \bty - \byhat)' \, \bSigma \, ( \bty - \byhat) \le \gamma^2 ,
\hspace{0.051cm} \Big\},
\end{equation}
where
\begin{equation}
    \bSigma = \sum_{j=1}^m \lambda_{j} v_j^Tv_j. 
\end{equation}
Let $l_j$ be the length of axis $j$ such that $\lambda_j = \frac{1}{l_j}$. 
Considering $\Gamma$ as the budget of uncertainty, the budgeted ellipsoidal uncertainty set can be derived as 
\begin{equation}
\label{ellipsoide U set}
\cU_e = \Big\{ \hspace{0.051cm} \bty \in \mathbb{R}^{{m}} : \hspace{0.2cm}  \bty = \byhat + \sum_{j=1}^m \mu_jv_j , \hspace{0.07cm}  \sum_{j=1}^m (\frac{\mu_j}{l_j})^2 \le \gamma^2 , \hspace{0.051cm}  \sum_{j=1}^m |\mu_j| \le \Gamma , \hspace{0.051cm}  \mu_j \in \mathbb{R}
\hspace{0.051cm} \Big\}.
\end{equation}
where the nominal vector $\byhat$ denotes the center of the ellipsoid, $\bty$ captures the deviation along the axes, and $\gamma \in [0,\gamma^{max}]$ is a controllable parameter that captures the maximum allowed  deviations in the uncertainty set.
The budget parameter $\Gamma \in [0,\Gamma^{max}]$, limits the total deviation of the ellipsoidal axes from their nominal values, where $\Gamma^{max}=\sum_{j=1}^m l_j$ corresponds to the full budget. 
Considering $\cU_e$ as the budgeted ellipsoidal uncertainty set, the robust problem with budget of uncertainty \ref{M2 with Budget} can be reformulated as \citep{kurtz2018robust}.
\begin{subequations}\label{ellipsoid with Budget}
\begin{align}
[\textbf{R($\gamma, \Gamma$)}]: \quad \min_{\bx,\by} \quad & \Big\{ \mathbf{c}_1'\bx+ \max_{\boldsymbol{\mu}} \mathbf{c}_2'\Big(
\byhat + \sum_{j=1}^m \mu_jv_j \Big) \Big\}, & \label{eM2-a budget} \\
\text{s.t} \quad & \bA\bx+\bB\by \le \bg, & \label{eM2-b budget}\\ 
& \by \le \byhat + \sum_{j=1}^m \mu_jv_j , & \label{eM2-d budget} \\
& \sum_{j=1}^m (\frac{\mu_j}{l_j})^2 \le \gamma^2,  &\label{em2- cone} \\ 
& \sum_{j=1}^m |\mu_j| \le \Gamma, & \\
& \boldsymbol{\mu} \in \mathbb{R}^m, &  \label{eM2-gamma}\\
& \bx,\by \ge \bzero. & \label{eM2-e budget}
\end{align}
\end{subequations}
Constraint \eqref{em2- cone} is a second-order cone defining the ellipsoidal feasible region for the uncertain parameters within radius $\gamma$. Similar to the polyhedral case, the budget of uncertainty may not be fully utilized, which may lead to an ineffective uncertainty set that results in a conservative solution.

The ellipsoidal uncertainty set $\cY_e$ is a nonempty bounded closed convex set, so the methodology proposed in Section~\ref{sec:extension} applies. We first define the axis-aligned minimum bounding box $\cH$ as 
\[\cH=\Big\{ \hspace{0.051cm} \bty \in \mathbb{R}^{{m}} : \quad  \buy \le \bty \le \boy, \quad \buy = \min_{\bty \in \cY_e}, \quad  \boy = \max_{\bty \in \cY_e} \Big\}.\]
We can then use the $\buy$ and $\boy$ to find the admissible interval $[\bus,\bos]$ using formulation~\eqref{sub}. Figure~\ref{fig: eff ellipsoid}(a) shows an example of an ineffective budget of uncertainty where the uncertain parameter $\tilde{y}_1$ cannot be utilized more than $\overline{s}_1$. Therefore, given the objective function shown with the purple line in Figure~\ref{fig: eff ellipsoid}, at point A, there is an ineffective allocation of the budget and increasing the budget 
would not impact the robust solution $y_1$. On the other hand, Figure~\ref{fig: eff ellipsoid}(b) shows an \emph{effective ellipsoid} (orange ellipsoid), which only considers effective uncertainty. In this case, the effective worst-case scenario under the same objective function corresponds to point B which leads to a less conservative solution 
since the effective worst-case scenario allocates the budget of uncertainty more in favor of uncertain parameter $\tilde{y}_2$, which can potentially be utilized more than $\tilde{y}_1$. Proposition~\ref{initial prop ellipsoid} formalizes the conditions that result in ineffective budgets in ellipsoidal uncertainty sets, and similarly, the possible cases of the admissible interval with respect to the ellipsoidal uncertainty are outlined in Remark~\ref{prop - ellipsoide max admissiblity}.

\begin{figure}
    \centering
    \includegraphics[width=0.9\textwidth, height=0.3\textheight]{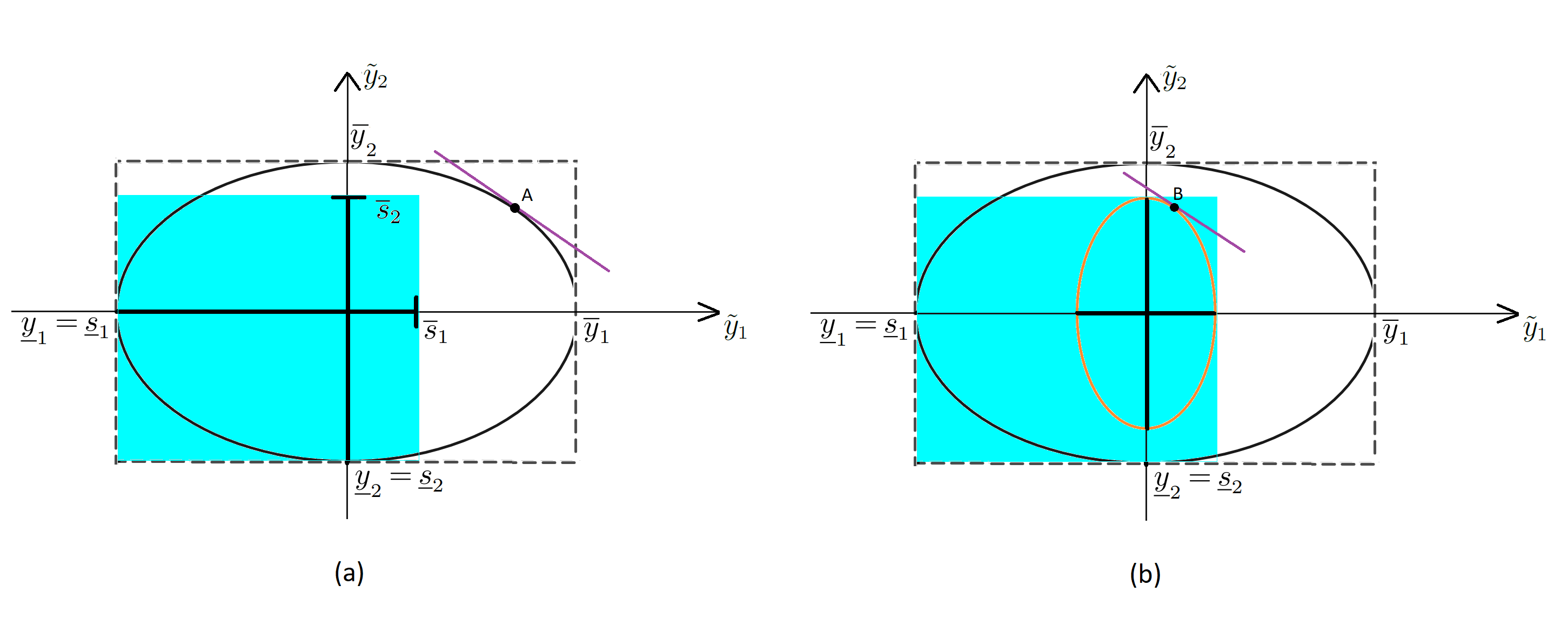}
    \caption{(a) Admissible area (shaded area) and (b) effective uncertainty set (orange ellipsoid) of the initial ellipsoidal uncertainty set. The purple line is the objective function of inner maximization problem in \eqref{eM2-a budget} and the dotted lines represent the axis-aligned minimum bounding box $\mathcal{H}$. 
    }
        \label{fig: eff ellipsoid}
\end{figure}

\begin{proposition}\label{initial prop ellipsoid}
For budgets of uncertainty $\Gammaz \geq \Gammal$ and radius $\gamma^0 \ge \gamma^l$, problems $[\mathbf{R}(\gamma^0,\Gamma^0)]$ and $[\mathbf{R}(\gamma^l,\Gamma^l)]$  with optimal solutions $(\bxz, \byz,\boldsymbol{\mu}^0)$ and $(\bxl, \byl, \boldsymbol{\mu}^0)$, respectively, result in the same optimal values $\bxz=\bxl$ and $\byz=\byl$ if \, 
$\exists \, i \in \{1, \dots, n\} \, \mid \, c_{2,(i)} \ge \max_k \hspace{0.1cm} c_{1,(k)}, \bA_i \bxl \geq 0$ and $\bA_i \bxl + \bB_i \byl = g_i$.    
\end{proposition}
\proof{Proof.}
First, for $\gamma^0 \ge \gamma^l$, 
the ellipsoid corresponding to $\gamma^0$ contains that of $\gamma^l$. Also, since $\Gamma^0 \ge \Gamma^l$, more deviations is allowed for $\Gamma^0$. Thus it must be that $\boldsymbol{\mu}^0 \geq \boldsymbol{\mu}^l$ and 
$\byz \geq \byl$. The rest of the proof follows from the proof of Proposition~\ref{initial prop}.
\Halmos
\endproof \vspace{1em}

\begin{remark}
\label{prop - ellipsoide max admissiblity}
For any $j$, the admissible interval $[\underline{s}_j,\overline{s_j}]$ can always be categorized as one of the following four cases:
\begin{enumerate}
\item $\suj=\yuj$ and \hspace{0.05cm} $ |a_j| < \soj \le \yoj$
\item $\suj=\yuj$ and \hspace{0.05cm} $\yhj < \soj \le  |a_j| $
\item $\suj=\yuj$ and \hspace{0.05cm} $\yuj < \soj \le \yhj$
\item $\suj=\soj \le \yuj$
\end{enumerate}
where $(\hat{y}_k, a_j=\hat{y}_j \pm (\gamma l_j |v_j|), k \in \{1,\dots,m\}/\{j\}$  show the intercepts of the ellipsoid with axis $j$.  
\end{remark}
\proof{Proof.}
Follows from Proposition \ref{prop - admissible interval}. \Halmos

Cases 1 and 2 of Remark~\ref{prop - ellipsoide max admissiblity} correspond to the ellipsoid axes that can be fully or partially utilized within radius $\gamma$. Cases~3 and~4 are identical to those of Proposition~\ref{prop - admissible interval} and correspond to the ineffective ellipsoid axes where the budget $\Gamma$ has no impact on them. 
Therefore, we can reformulate the budgeted-ellipsoidal uncertainty set by considering the budget of uncertainty $\Gamma$ only for the effective deviations corresponding to cases 1 and 2 and exclude the ineffective deviations corresponding to cases 3 and 4. To ensure the feasibility of case 3 and case 4, we impose constraint $\bty \le \bos$, without assigning any budget to them. In the effective ellipsoid, the length of each ellipsoid direction $v_j$ is recalculated based on its maximum admissibility. For vector $j$, $l^e_j = (\overline{s}_j-\hat{y}_j)$ is the effective length of vector $v_j$, which will widen the ellipsoid in the direction of axes with higher admissibility.
\begin{definition} \label{Def2} For $j\in M=\{1,\dots, m\}$, let $I := \{j : \overline{s}_j > \hat{y}_j\}$ and $I^c := M/I$. Dividing vectors based on $I$ and $I^c$, let $\bts=(\bts^{I},\bts^{I^c})$, $\bos=(\bos^{I},\bos^{I^c})$, and $\byhat=(\byhat^{I},\byhat^{I^c})$. The effective budgeted-ellipsoidal uncertainty  $\cU^{E}_e$ is defined as
\begin{equation}
\label{effective ellipsoide U set}
\cU^{E}_e = \Big\{ \hspace{0.051cm} \bts^{I} \in \mathbb{R}^{{|I|}}, \bts^{I^c} \in \mathbb{R}^{{m-|I|}} : \hspace{0.2cm}  \bts^I = \byhat^I + \sum_{j \in I} \mu_jv_j, \hspace{0.07cm}  \sum_{j \in I} (\frac{\mu_j}{l^e_j})^2 \le \gamma^2 , \hspace{0.051cm}  \sum_{j \in I} |\mu_j| \le \Gamma , \hspace{0.051cm}  \mu_j \in \mathbb{R}
\hspace{0.051cm}, \hspace{0.051cm}  \bts^{I^c} \le \bos^{I^c},  \Big\},
\end{equation}
where $\forall j \in I$, parameter $l^e_j=(\soj - \hat{y}_j)$ is the length of effective vector $v_j$.
\end{definition} 

\begin{proposition}
The effective budgeted-ellipsoidal uncertainty set $\cU^{E}_e$ eliminates the ineffective part of the uncertainty set $\cU_e$.  
\end{proposition}
\proof{Proof.}
Entries of the set $I$ represent the effective vectors $v_j$. The first four terms of $\cU^{E}_e$ define a budgeted-ellipsoid only using the effective axes. Also, for $j\in I^c$, the term $\bts^{I^c} \le \bos^{I^c}$ ensures the feasibility of the problem for the ineffective axes.\Halmos
\vspace{0.3cm}

Finally, given sets $j\in M=\{1,\dots, m\}$, $I := \{j : \overline{s}_j > \hat{y}_j\}$ and $I^c := M/I$ from the first stage, we rewrite vectors $\byhat=(\byhat^I,\byhat^{I^c}), \by=(\by^I,\by^{I^c}), \bts=(\bts^I,\bts^{I^c})$, and $\bos=(\bos^I,\bos^{I^c})$. Substituting  $\bts^I,\bts^{I^c}$, the second-stage problem is
\begin{subequations}\label{second stage ellipsoid with Budget}
\begin{align}
\min_{\bx,\by} \quad & \Big\{ \mathbf{c}_1'\bx+ \max_{\boldsymbol{\mu}} \mathbf{c}_2'\Big(
\byhat^I + \sum_{j \in  I} \mu_jv_j \Big) \Big\}, & \label{seM2-a budget} \\
\text{s.t} \quad & \bA\bx+\bB\by \le \bg, & \label{seM2-b budget}\\ & \by^{I^c} \le \bos^{I^c} \\
& \by^{I} \le \byhat^I + \sum_{j \in  I} \mu_jv_j , & \label{seM2-d budget} \\
& \sum_{j \in I} (\frac{\mu_j}{l_j})^2 \le \gamma^2 \label{sem2- cone} \\ & \sum_{j \in I} |\mu_j| \le \Gamma, \\& \boldsymbol{\mu} \in \mathbb{R}^m
  \label{seM2-gamma}
\\ & \bx,\by \ge \bzero. &
\label{seM2-e budget}
\end{align}
\end{subequations}

The second-stage is not linear and its robust counterpart optimization formulation is equivalent to a second-order cone program that can be reformulated using existing methods in the literature \citep{alizadeh2003second,bertsimas2011theory,kurtz2018robust}. 
}

\section{The SCED Problem} \label{appendix SCED}
In this section, we present three variants of modeling the multi-period security-constrained economic dispatch (SCED) problem based on settings from \cite{li2016adaptive}. In the SCED models presented here, parameters $\underline{W}_{k,t}$, $\hat{W}_{k,t}$, $\tilde{W}_{k,t}$, and $\overline{W}_{k,t}$ corresponds to the uncertain vectors $\buy$, $\byhat$, $\bty$ and $\boy$ used in this paper, respectively. In what follows, we first present all notations used in the SCED models. Next, we present the three SCED optimization models used in the numerical results, namely, the deterministic model, the conventional robust model, and the proposed two-stage robust model. 

\begingroup\makeatletter\def\f@size{10}\check@mathfonts 
\printnomenclature[1.6cm]
\nomenclature[A,04]{$\mathcal{G}_i$}{Set of conventional generators at bus \textit{i}}
\nomenclature[A,05]{$\mathcal{K}_i$}{Set of wind farms at bus \textit{i}}
\nomenclature[A,06]{$\mathcal{F}$}{Set of transmission lines}
\nomenclature[a,02]{$\mathcal{N}$}{Set of buses in power system}
\nomenclature[a,07]{$\mathcal{T}$}{Set of scheduling periods}
\nomenclature[B,03]{$D_{i,t}$}{Load demand at bus \textit{i} during period \textit{t}}
\nomenclature[B,08]{$\underline{P}_{g}/\overline{P}_{g}$}{Upper/Lower output limit of conventional generator \textit{g}}
\nomenclature[B,09]{$\hat{W}_{k,t}$/$\tilde{W}_{k,t}$}{Nominal/Uncertain parameter of available wind power output of wind farm \textit{k} during
period \textit{t}}
\nomenclature[B,04]{$\overline{F}_{f}/\underline{F}_{f}$}{Upper/Lower power flow limit of transmission line \textit{f}}
\nomenclature[B,06]{$G_{f,i}$}{Sensitivity matrix of power injection at internal bus \textit{i} to flow on line \textit{f}}
\nomenclature[B,13]{${R}^{u}_{t}$/${R}^{d}_{t}$}{Upward/downward spinning reserve requirement during period \textit{t}}
\nomenclature[B,12]{${U}_{g}$/${D}_{g}$}{Upward/downward ramping rate of conventional generator \textit{g}}
\nomenclature[B,02]{$\sigma_k$}{Penalty cost of wind power curtailment of wind farm \textit{k}
}
\nomenclature[B,01]{$C_{g}$}{Generation cost of conventional generator \textit{g}}
\nomenclature[E,06]{$r^+_{g,t}$/$r^-_{g,t}$}{Upward/downward spinning reserve of conventional generator \textit{g} during period \textit{t}}
\nomenclature[E,01]{$p_{g,t}$}{Power output of conventional generator \textit{g} during period \textit{t}}
\nomenclature[E,05]{$p^{W}_{k,t}$}{Power output of wind farm \textit{k} during period \textit{t}}

\noindent {\bf The Deterministic Model}

\noindent The deterministic model considers the predicted wind power and assumes the prediction is perfect. It then finds the power dispatches based on the predicted wind power output, since the deterministic model does not account for wind uncertainty. The mathematical formulation of the multi-period deterministic SCED problem is as follows:
\begin{align}
\min_{
} \quad & \sum_{t \in \mathcal{T}}\sum_{i \in \mathcal{N}} \left( \sum_{g \in \mathcal{G}_i} C_gp_{g,t} + \sum_{k \in \mathcal{K}_i} \sigma_k \Big( \hat{W}_{k,t}-p^W_{k,t} \Big) \right), & \label{SCED: OBJ} \\
\text{s.t:} \quad & \sum_{i \in \mathcal{N}}\Big( \sum_{g \in \mathcal{G}_i} p_{g,t} + \sum_{k \in \mathcal{K}_i} p^W_{k,t} \Big) = \sum_{i \in \mathcal{N}} D_{i,t} , & \forall t \in \mathcal{T}
\label{AC grid: power balance} \\
&\underline{F}_{f} \leq \sum_{i \in \mathcal{N}} G_{f,i} \left( \sum_{g \in \mathcal{G}_i} p_{g,t} + \sum_{k \in \mathcal{K}_i} p^W_{k,t} - D_{i,t} \right) \leq \overline{F}_{f}, & \forall f \in \mathcal{F}, \forall t \in \mathcal{T}\label{AC grid:internal line limits} 
\\
&\sum_{i \in \mathcal{N}} \sum_{g \in \mathcal{G}_i} r^{+}_{g,t} \ge R^{u}_{t} , & \forall t \in \mathcal{T} \label{remaining positive spinning reserve in the system} \\
&\sum_{i \in \mathcal{N}} \sum_{g \in \mathcal{G}_i} r^{-}_{g,t} \ge R^{d}_{t} ,& \forall t \in \mathcal{T}\label{remaining negative spinning reserve in the system} \\
&0\le r^{+}_{g,t} \le \min \left\lbrace \overline{P}_{g} - p_{g,t}\hspace{0.05cm}, \hspace{0.05cm} U_{g}.\Delta t  \right\rbrace,&\forall g \in \mathcal{G}_i, i \in \mathcal{N}, t \in \mathcal{T}
\label{AC grid: final spinning up}\\
&0\le r^{-}_{g,t} \le \min \left\lbrace p_{g,t} - \underline{P}_{g}\hspace{0.05cm},\hspace{0.05cm} D_{g}.\Delta t \right\rbrace, & \forall g \in \mathcal{G}_i, i \in \mathcal{N}, t \in \mathcal{T}
\label{AC grid: final spinning down} \\
&-U^d_g.\Delta t \le p_{g,t} - p_{g,t-1} \le U^u_g.\Delta t, & \forall g \in \mathcal{G}_i, i \in \mathcal{N}, t \in \mathcal{T}
 \label{AC grid:ramping rate}
\\
&\underline{P}_{g} \leq p_{g,t} \leq \overline{P}_{g}, & \forall g \in \mathcal{G}_i, i \in \mathcal{N}, t \in \mathcal{T}
\label{AC grid: PG limits}
\\
& 0 \le p^W_{k,t} \leq \hat{W}_{k,t}. & \forall k \in \mathcal{K}_i, i \in \mathcal{N}. t \in \mathcal{T} \label{AC grid: PW limits}
\end{align}
\endgroup

The objective function \eqref{SCED: OBJ} minimizes the total operational cost consisting of generation cost of conventional generators and wind power curtailment cost. Constraint \eqref{AC grid: power balance} ensures the power balance between generation and load. In constraint \eqref{AC grid:internal line limits}, the power flows limits are considered. The requirements for upward and downward spinning reserve are considered in constraints \eqref{remaining positive spinning reserve in the system} and \eqref{remaining negative spinning reserve in the system}, respectively.  Constraints \eqref{AC grid: final spinning up} and \eqref{AC grid: final spinning down} address the capacity of upward and downward reserves of each conventional generator. Constraint \eqref{AC grid:ramping rate} ensures that the dynamic changes in power outputs are limited with respect to the ramping rate of conventional generators. The generation limits of conventional generators and wind farms are shown in \eqref{AC grid: PG limits} and \eqref{AC grid: PW limits}, respectively. 

\noindent {\textbf{The Conventional Robust Model}}

\noindent The conventional robust model accounts for uncertainty but does not account for effective budgets. The objective function of the conventional robust approach is as follows \citep{li2016adaptive}
\begin{equation}\label{SCED: OBJ robust conventional}
\hspace{-0.7cm}  \min_{
} \hspace{0.3cm} \sum_{t \in \mathcal{T}}\sum_{i \in \mathcal{N}} \left(  \sum_{g \in \mathcal{G}_i} C_gp_{g,t} + \max_{\tilde{W}_{k,t}} \sum_{k \in \mathcal{K}_i} \sigma_k \Big( \tilde{W}_{k,t}-p^W_{k,t} \Big) \right),
\end{equation}
where the inner maximization finds the wind power curtailment under the absolute worst-case scenario, and then the optimal solution is obtained under this worst-case scenario.
In the conventional robust model, constraint \eqref{AC grid: PW limits con. robust} corresponds to the wind power dispatch limit, and constraints \eqref{USet1 Wtilde}-\eqref{Uset3 scaled deviations} correspond to a polyhedral uncertainty set incorporating the budget of uncertainty $\Gamma_t$ in the model.
\begin{align}
&0 \le p^W_{k,t} \leq \tilde{W}_{k,t},  & \forall k \in \mathcal{K}_i, i \in \mathcal{N}, t \in \mathcal{T}
\label{AC grid: PW limits con. robust}
\\
&\tilde{W}_{k,t} = \hat{W}_{k,t} + z^+_{k,t}(\overline{W}_{k,t}-\hat{W}_{k,t})+z^-(\underline{W}_{k,t}-\hat{W}_{k,t}), 
& \forall k \in \mathcal{K}_i, i \in \mathcal{N}^{I}, t \in \mathcal{T} \label{USet1 Wtilde} \\
&\displaystyle \sum_{i \in \mathcal{N}^{I}_a}  \sum_{k \in \mathcal{K}_i} (z^-_{k,t}+z^+_{k,t}) \le \Gamma_{t},& \forall t \in \mathcal{T} \label{Uset2 total deviations} \\
&0\le z^+_{k,t} \le 1, & \forall k \in \mathcal{K}_i, i \in \mathcal{N}^{I}, t \in \mathcal{T} \label{Uset3 scaled deviations} \\
&0\le z^-_{k,t} \le 1. & \forall k \in \mathcal{K}_i, i \in \mathcal{N}^{I}, t \in \mathcal{T} \label{Uset4 scaled deviations} 
\end{align}

Using Duality Theorems, a linear reformulation of the conventional roust model is obtained as follows where $\boldsymbol{\xi}$ and $\boldsymbol{\mu}$ are dual vectors corresponding to  \eqref{Uset2 total deviations} and \eqref{Uset3 scaled deviations}. Note that the constraint corresponding to \eqref{Uset4 scaled deviations} is removed as it becomes redundant, shown in Proposition \ref{prop worst-case realization interval}:
\begin{align}
\hspace{-0.3cm} \min_{
}
\quad & \sum_{t \in \mathcal{T}}\sum_{i \in \mathcal{N}} \left( \sum_{g \in \mathcal{G}_i}  C_{g}p_{g,t} +
\sum_{k \in \mathcal{K}_i} \sigma_k \big( \hat{W}_{k,t} + \mu_{k,t} + \Gamma_{t} \xi_{t}-{p}^{W}_{k,t} \big) \right), & \label{conv robust - dual obj of Stage I}\\
\text{s.t:} \quad &\mu_{k,t} + \xi_{t} \ge (\overline{W}_{k,t}-\hat{W}_{k,t}), &\forall k \in \mathcal{K}_i, i \in \mathcal{N}^I, t \in \mathcal{T} \\
&\boldsymbol{\mu}, \boldsymbol{\xi} \ge \mathbf{0}\\
& \eqref{AC grid: power balance}-\eqref{AC grid: PG limits}, \eqref{AC grid: PW limits con. robust}-\eqref{Uset4 scaled deviations}.
\nonumber
\end{align}
\noindent {\bf The Proposed Robust Model}

\noindent The proposed robust model considers effective budgets of uncertainty and uses the proposed two-stage approach. 

\noindent
{ \bf Stage (I):}
The auxiliary optimization problem \eqref{stage I SCED} is used to find the largest admissible wind power interval $[\underline{s}_{k,t},\overline{s}_{k,t}]$:
\allowdisplaybreaks
\begin{subequations} \label{stage I SCED}
\begin{align} 
\min_{
}  \quad & \sum_{t \in \mathcal{T}} \sum_{i \in \mathcal{N}^I}\sum_{k \in \mathcal{K}_i} \big( \overline{W}_{k,t}-\overline{s}_{k,t}\big) + \big(\underline{W}_{k,t} - \underline{s}_{k,t}\big) & \label{stage I SCED - Obj} \\
\text{s.t.} \quad & \sum_{i \in \mathcal{N}^{I}} G^{I}_{f,i} \left( \sum_{g \in \mathcal{G}_i} p_{g,t} - D_{i,t} + \sum_{k \in \mathcal{K}_i} \underline{s}_{k,t}\right) + \sum_{i \in \mathcal{N}^{I}}
\sum_{k \in \mathcal{K}_i} \alpha_{k,f,t} \leq \overline{F}_{f}, & \label{Robust counterpart 1 for upper limit of the power flows}\\
&  \hspace{7.3cm} \forall k \in \mathcal{K}_i, i \in \mathcal{N}^I, f \in \mathcal{F}, t \in \mathcal{T} \nonumber\\
& \alpha_{k,f,t} \geq G^I_{f,i}(\overline{s}_{k,t} -\underline{s}_{k,t}), \qquad \quad \qquad  \hspace{1.55cm} \forall k \in \mathcal{K}_i, i \in \mathcal{N}^I, f \in \mathcal{F}, t \in \mathcal{T} \label{Robust counterpart 2 for upper limit of the power flows} \\
& \sum_{i \in \mathcal{N}^{I}} G^{I}_{f,i} \left( \sum_{g \in \mathcal{G}_i} p_{g,t} - D_{i,t} + \sum_{k \in K_i} \underline{s}_{k,t}\right) - \sum_{i \in \mathcal{N}^{I}}
\sum_{k \in \mathcal{K}_i} \zeta_{k,f,t} \geq \underline{F}_{f}, \label{Robust counterpart 1 for lower limit of the power flow constraint} \\
&  \hspace{7.3cm} \forall k \in \mathcal{K}_i, i \in \mathcal{N}^I, f \in \mathcal{F}, t \in \mathcal{T} \nonumber\\
& \zeta_{k,f,t} \geq -G^I_{f,i}(\overline{s}_{k,t} -\underline{s}_{k,t}), \hspace{2.5cm} \quad \quad \forall k \in \mathcal{K}_i, i \in \mathcal{N}^I, f \in \mathcal{F}, t \in \mathcal{T}\label{Robust counterpart 2 for lower limit of the power flow constraint} \\
& \sum_{i \in \mathcal{N}^{I}} \left( \sum_{g \in \mathcal{G}_i} p_{g,t} + \sum_{g \in \mathcal{G}_i} r^{+}_{g,t} - D_{i,t} + \sum_{k \in \mathcal{K}_i} \underline{s}_{k,t} - \sum_{k \in \mathcal{K}_i} \eta_{k,t} \right) \geq R^{u}_{t},    \hspace{1.0cm} \forall t \in \mathcal{T}  \label{Robust counterpart 1 for positive spinning reserve constraint} \\
& \eta_{k,t} \geq -(\overline{s}_{k,t} -\underline{s}_{k,t}), \hspace{4.4cm} \qquad \forall k \in \mathcal{K}_i, i \in \mathcal{N}^I, t \in \mathcal{T} \label{Robust counterpart 2 for positive spinning reserve constraint} \\
& \sum_{i \in \mathcal{N}^{I}} \left(  \sum_{g \in \mathcal{G}_i} p_{g,t} - \sum_{g \in \mathcal{G}_i} r^{-}_{g,t} - D_{i,t} + \sum_{k \in \mathcal{K}_i} \underline{s}_{k,t} + \sum_{k \in \mathcal{K}_i} \beta_{k,t} \right) \le R^{d}_{t}, \hspace{1.0cm} \forall t \in \mathcal{T} 
\label{Robust counterpart 1 for negative spinning reserve constraint} \\
& \beta_{k,t} \geq \overline{s}_{k,t} -\underline{s}_{k,t}, \hspace{5cm} \qquad \forall k \in \mathcal{K}_i, i \in \mathcal{N}^I, t \in \mathcal{T}  \label{Robust counterpart 2 for negative spinning reserve constraint} \\
& \overline{s}_{k,t} \le \overline{W}_{k,t}, \hspace{5.76cm} \qquad \forall k \in \mathcal{K}_i, i \in \mathcal{N}^{I}, t \in \mathcal{T} \label{new initial relationships of Wmax Wmin and Smax Smin1} \\
& \underline{s}_{k,t} \le \underline{W}_{k,t}, \hspace{6.5cm} \forall k \in \mathcal{K}_i, i \in \mathcal{N}^{I}, t \in \mathcal{T} \label{new initial relationships of Wmax Wmin and Smax Smin2}\\
& \bm {\alpha, \zeta, \eta, \beta, \mathbf{r^+}, \mathbf{r^-}, \bos, \bus, \mathbf{p}
} \ge \bm{0}, \label{SCED dual variables}
\end{align}
\end{subequations}
where vectors $\bm {\alpha,  \zeta,\eta, \beta
}$ are auxiliary variables used in the robust counterpart reformulation. 
Given the admissible interval $[\underline{s}_{k,t},\overline{s}_{k,t}]$, the interval $[\hat{s}_{k,t},\overline{s}_{k,t}]$ is obtained and used in the uncertainty set to represent the effective worst-case scenario.

\noindent
{ \bf Stage (II):} Given the optimal solution of Stage I, the effective budget of uncertainty $\Gamma^{E}$
is obtained and incorporated in the Stage II formulation as follows:
\begin{align}\label{dual obj of Stage II} \min_{
} \quad & \left\{ \sum_{t \in \mathcal{T}} \sum_{i \in \mathcal{N}^I} \left( \sum_{g \in \mathcal{G}_i}  C_{g}p_{g,t} + \sum_{k \in \mathcal{K}_i} \sigma_k \big( \hat{s}_{k,t} + \mu_{k,t} -v_{k,t}\lambda_{k,t} + \Gamma^{E}_{t} \xi_{t}-{p}^{W}_{k,t} \big) \right)
\right\}&\\
\text{s.t.} \quad &\mu_{k,t} -\lambda_{k,t} + e_{k,t}\xi_{t} \ge (\overline{s}_{k,t}-\hat{s}_{k,t}), \hspace{1.5cm} \qquad \qquad \forall k \in \mathcal{K}_i, i \in \mathcal{N}^I, t \in \mathcal{T} \label{dual - obj const II} \\
&p^W_{k,t} \leq \tilde{s}_{k,t},\\
& \tilde{s}_{k,t} = \hat{s}_{k,t} + r_{k,t}(\overline{s}_{k,t}-\hat{s}_{k,t}), \hspace{3.8cm} \forall k \in \mathcal{K}_i, i \in \mathcal{N}^{I}, t \in \mathcal{T} \\
& \sum_{i \in \mathcal{N}^{I}}  \sum_{k \in \mathcal{K}_i} e_{k,t}.r_{k,t} \le \Gamma^{E}_{t}, \hspace{6.8cm} \forall t \in \mathcal{T} \label{SCED Def 1 - b}\\ & v_{k,t} \le r_{k,t} \le 1,  \hspace{8.24cm} \forall t \in \mathcal{T} \label{SCED Def 1 - c} \\
&\boldsymbol{\mu}, \boldsymbol{\lambda}, \boldsymbol{\xi} \ge \mathbf{0}, & \quad \label{non-negative dual - obj const II}\\
& \eqref{AC grid: power balance}-\eqref{AC grid: PG limits}
\end{align}
where parameters $e_{k,t}$, $v_{k,t}$, and $\Gamma^{E}_{t}$ are calculated based on Section \ref{ch2 effective budget}.
%

{
\section{Examples of other Applications} \label{appendix other applications}
In this section, we provide two small numerical instances of the applications mentioned in Section~\ref{sec other applications}. 
\subsection{Patient Scheduling} 
\label{appendix patient scheduling}
Consider a healthcare system with multiple patient priorities $p \in \mathcal{P}=\{1,\dots,P\}$ 
and let parameters $s_p$ and $w_p$ denote the service time and the maximum wait-time target for each priority. Let $g_{t,p}$ be the penalty cost for patients with priority level $p$ who arrive on day $t$ and do not receive service within the time horizon. 
Let $\tau=\{1,\dots,T\}$ 
be the set of time periods, 
Let $\tilde{d}_{t,p} \in [\,\underline{d}_{t,p},\, \overline{d}_{t,p}\,]$ be the uncertain number of patients of priority level $p$ that arrive on day $t$. 
The decision variables in this problem are $x_{t,n,p}$ and $c_{i,t}$ which correspond to the number of patients with priority level $p$ who receive service on day $t$ after waiting for $n$ days, and 
capacity on day $t$, respectively. The capacity has an associated piece-wise linear cost of $\theta_i$ for each linear segment $i$ with a given limit of $L_i$. We assume that $\theta_1=0$, meaning that there is no cost for using the existing capacity of $L_1$. 
Considering binaries $q_{i,t} \in \{0,1\}$, the mathematical formulation of this patient and capacity scheduling problem 
under demand uncertainty is as follows \citep{Shahraki2020advance}:

\begin{subequations}\label{Inv 1}
\begin{align}
\quad \min_{\mathbf{q,x,c}}  \, & \Bigg\{ \sum_{t =1}^T \sum_{i =1}^K \theta_i c_{i,t} + 
\max_{\mathbf{z}} \left\{\sum_{p=1}^{P}\sum_{t=1}^{T} \Big( \hat{d}_{t,p} + z_{t,p}(\overline{d}_{t,p}-\hat{d}_{t,p}) - \sum_{i=0}^{T-t} x_{T-i,T-t-i+1,p}   \Big)g_{t,p} \right\} \Bigg\},  \label{PS obj} \\
\text{s.t} \quad 
& \sum_{p=1}^{P}\sum_{n=1}^{t} s_px_{t,n,p} \le \sum_{i=1}^K c_{i,t}, \qquad \qquad \qquad \qquad \forall t \in \tau  \label{PS 2}\\
& (L_i-L_{i-1})q_{i,t} \le c_{i,t} \le (L_i-L_{i-1})q_{i-1,t} \qquad \forall t \in \tau, \forall i\in \{2, \dots, K-1\}   \label{PS 4}\\
& 0 \le c_{i,t} \le L_i \hspace{4.86cm} \forall t \in \tau , \forall i \in \mathcal{K} \label{PS 5}\\
& \sum_{i=0}^{n-1} x_{t-i,n-i,p} \le \hat{d}_{t,p} + z_{t,p}(\overline{d}_{t,p}-\hat{d}_{t,p}), \qquad \hspace{0.3cm} \forall n \in \{1,...,t\}, \forall t \in \tau, \forall p \in \mathcal{P}  \label{PS 1}\\ 
&\sum_{t=1}^{T} z_{t,p} \le \Gamma_p, \hspace{4.6cm}  \forall p \in \mathcal{P} \label{PS -gamma} \\ 
& \bzero \le \mathbf{z} \le \bone,  \label{ps f z limits} \\ 
& \mathbf{x} \ge \bzero, \quad \text{ and integer}, \quad \mathbf{q} \in \{0,1\}^{K\tau}. 
\label{ps -e budget}
\end{align}
\end{subequations}
The objective function~\eqref{PS obj} optimizes the total cost, i.e., the cost of increasing capacity as well as the penalty cost of unserved patients, under the worst-case scenario. Constraint~\eqref{PS 2} restricts the total number of patients who receive service each day within the level of capacity for that day. Constraints~\eqref{PS 4}--\eqref{PS 5} determine the level of capacity for each day. Constraint~\eqref{PS 1} limits the number of patients receiving service to the demand. Finally, constraint~\eqref{PS -gamma} sets the budget of uncertainty. The readers are referred to~\cite{mahmoudzadeh2020robust} and~\cite{Shahraki2020advance} for further details on this problem. Grouping constraints \eqref{PS 2}-\eqref{PS 5}, this model can be reduced to formulation~\eqref{M2 with Budget}.

Table~\ref{Table data patient sch} shows the data used in our numerical example as well as the optimal solution using both the conventional and the proposed robust budget-of-uncertainty approaches.

\begin{table}[htbp]

\centering
\footnotesize
\caption{Data and the optimal solution of the numerical example on patient scheduling}
\label{Table data patient sch}
\begin{tabular}{cc|ccc}
\hline 
\multicolumn{2}{c|}{Data} & \multicolumn{3}{c}{Optimal Solution} \\ \hline
Parameters & Settings & Variables & Conventional approach & Proposed approach \\ \hline
$(P,T,K,s, \Gamma)$ & $(1,2,2,1,0.5)$ & $(z_{1,1},z_{2,1})$ & (0,0.5) & (0.5,0) \\ $\boldsymbol{\theta}$ & $(0,25)$ & $(q_{1,1},q_{1,2},q_{2,1},q_{2,2})$ & (1,1,0,0) & (1,1,0,0) \\ $\mathbf{g}$&$(2,2)$ & $(c_{1,1},c_{1,2},c_{2,1},c_{2,2})$ & (15,15,0,0) & (15,15,0,0) \\
$\mathbf{CL}$ & (15,35,50)  & $(x_{1,1,1},x_{2,1,1},x_{2,2,1})$ & (10,15,0) & (15,15,0) \\
$\hat{\mathbf{d}}$&$(10,30)$ & cost &25 &15 \\
$\overline{\mathbf{d}}$&(20,60)  \\
\hline
\end{tabular}
\end{table}%
In this numerical example, the conventional approach assigns the budget to the second day (i.e., $z_{2,1}=0.5$) which corresponds to the worst-case cost. However, this allocation may not be entirely utilized due to constraint \eqref{PS 2}. At optimality, we have $x_{2,1,1} + x_{2,2,1} \le 15$, meaning that the maximum demand that can be met on day 2 is 15. Since $\hat{d}_2=30$, even for a zero budget of uncertainty, $x_{2,1,1} + x_{2,2,1} \le 15 \le \hat{d}_2$ holds and thus assigning any budget of uncertainty to day 2 is ineffective. On the other hand, the proposed approach effectively reassigns the budget of uncertainty to day 1 at optimality, where the systems limits are sufficient to potentially satisfy the effective worst-case demand. Thus, this re-allocation of budget in favor of the effective worst-case scenario yields a less conservative solution and therefore corresponds to a lower cost, as shown in Table~\ref{Table data patient sch}.

\subsection{Inventory Planning}
\label{appendix inventory}
For this application, we use a basic version of the inventory model of \cite{bertsekas1995dynamic} and incorporate a RHS budget of uncertainty on the customer demand parameter. Let $I_k$ be the inventory level at the beginning of period $k \in K = \{1, ..., T\}$, $u_k$ the stock ordered at the beginning of period $k$, $v_k$ the binary variable indicating if there is an order in period k, and $\tilde{w}_k$ the uncertain demand in period $k$ which belongs to uncertainty set $[\underline{w}_k,\overline{w}_k]$. 
Let parameters $c^p_k, c^f_k, c^{sh}_k$, and $c^h_k$ denote the purchase, fixed, shortage, and holding cost, respectively. Also, let $I^{min}_k$ and $I^{max}_k$ be the lower and upper bound on inventory level $I_k$. The reformulated inventory planning model is 
\begin{subequations}\label{Inv 1}
\begin{align}
 \quad \min_{\mathbf{I,u,v}}  \qquad  & \Bigg\{ \sum_{k \in K}( c^p_ku_k + c^f_k v_k) + \max_{\bz} \Big\{ \sum_{k \in K} c^{sh}_k\big(\hat{w}_k + z_k(\overline{w}_k-\hat{w}_k)-I_k-u_k\big)\Big\} \nonumber \\ 
 & \qquad \hspace{2.2cm} + \max_{\bz} \Big\{\sum_{k \in K}c^h_k\Big(I_k+u_k-\hat{w}_k - z_k(\overline{w}_k-\hat{w}_k) \big)\Big\} \Bigg\}, \label{inv obj}\\
\text{s.t.} \qquad & u_k \le u_k^{max} v_k, & \forall k=1,...,\tau \label{Inv b budget}\\ 
& I_{k+1}-I_k-u_k+w_k=0, & \forall k=1,...,\tau \label{Inv bb budget} \\
& I^{min}_k \le I_k \le I_k^{max}, & \forall k=1,...,\tau+1 \label{Inv bbb budget} \\
& w_k \le \hat{w}_k + z_k(\overline{w}_k-\hat{w}_k), & \forall k=1,...,\tau\label{Inv d budget} \\
&\sum_{k=1}^{T} z_k \le \Gamma, &  \label{Inv -gamma} \\ 
& \bzero \le \mathbf{z} \le \bone, & \label{Inv f z limits}
\\ & \mathbf{u} \ge \bzero, \mathbf{v} \in \{0,1\}^\tau. &
\label{Inv -e budget}
\end{align}
\end{subequations}
In the objective function~\eqref{inv obj}, the first term corresponds to the fixed and purchase cost, and the inner inner maximization problems correspond to worst-case shortage and excess (holding) costs, respectively. Constraint~\eqref{Inv b budget} corresponds to purchase limit in each period. Constraints~\eqref{Inv bb budget} and~\eqref{Inv bbb budget} represent the inventory balance and inventory limit constraints, respectively. Constraint~\eqref{Inv d budget} is a reformulation of the demand constraint with uncertainty. In this problem, constraints~\eqref{Inv b budget} to~\eqref{Inv bbb budget} correspond to~\eqref{M2-b budget}, and constraint~\eqref{Inv d budget} corresponds to~\eqref{M2-d budget}. Hence, the inventory planning model~\eqref{Inv 1} reduces to formulation~\eqref{M2 with Budget} where $\max_{\bz} \mathbf{c}_2'\big(
\byhat + \bzn \odot (\buy - \byhat) + \bzp \odot (\boy - \byhat)-\by\big)$ in  formulation \eqref{M2 with Budget} corresponds to the worst-case shortage cost of problem \eqref{Inv 1}. 

Consider a small numerical example where the holding cost is high relative to the other costs. Table~\ref{Table data inventory} shows the data used in this numerical example as well as the optimal solution using both the conventional and proposed robust budget-of-uncertainty approaches.
\begin{table}[htbp]

\centering
\footnotesize
\caption{ Data and the optimal solution of the numerical example in inventory planning}
\label{Table data inventory}
\begin{tabular}{cc|ccc}
\hline 
\multicolumn{2}{c|}{Data} & \multicolumn{3}{c}{Optimal Solution} \\ \hline
Parameters & Settings & Variables & Conventional approach & Proposed approach \\ \hline
$(\tau,I_1,I_4,\Gamma)$ & $(3,0,0,2)$ & $(z_{1},z_{2},z_3)$ & (1,1,0) & (0.35,0,1)* \\ $\mathbf{c}^p$ & $(2,2,3)$ & $(v_{1},v_{2},v_3)$ & (1,1,1) & (1,1,1) \\ $\mathbf{c}^f$&$(1,1,1)$ & $(u_{1},u_{2},u_3)$ & (5,5,4) & (5,5,6) \\
$\mathbf{c}^{sh}$ & (4,4,4)  & $(I_1,I_2,I_3,I_4)$ & (0,-2,-2,0) & (0,-2,-2,0) \\$\mathbf{c}^{h}$ & (10,10,10) & ($w_1,w_2,w_3$) & (7,5,2) & (7,5,6) \\
$\mathbf{u}^{max}$&$(5,5,6)$ & cost &83 &41 \\
$\mathbf{I}^{min}$&$(-2,-2,-2)$  \\
$\mathbf{I}^{max}$&$(2,2,2)$ \\
$\hat{\mathbf{w}}$&$(5,5,2)$  \\
$\overline{\mathbf{w}}$&(10,8,4) & \multicolumn{3}{l}{*The new deviations are re-sclaed based on length of $\bos-\bshat$} \\
\hline 
\end{tabular}
\end{table}
%
For this dataset, the conventional approach finds the worst-case scenario by assigning the budget to a full deviation to period 1 and 2, and no deviation to period 3. However, due to limitations enforced by other constraints, some of these deviations are ineffective. For example, in period~1, it can be seen that $I_1+u_1 \le 5$ and thus $w_1\le 7$, meaning that  we cannot meet more than 7 units of demand, even though the full deviation allocated to period~1 assumes that the demand can go up to 10 units. The allocation of an extra 3 units of deviation to the demand of that period would not change the optimal solution, and is hence, ineffective, and would lead to an over-conservative solution. 
On the contrary, the proposed approach takes this ineffective deviation into account and 
reallocates the budget of uncertainty in favor of the effective worst-case scenario. In this case, the model allows for a full deviation of the uncertain parameter in period 3. Subsequently, the proposed approach leads to a less conservative solution and a lower cost, as shown in Table \ref{Table data inventory}.
%

We observed that when all holding costs are strictly less than the shortage costs, there was no shortage at optimality, and therefore, our proposed approach and the conventional budget-of-uncertainty approach both corresponded to the same solution, simply because there was no ineffective budget of uncertainty and both approaches could effectively control the solution conservatism in this particular case. Otherwise, there is always a potential for having an ineffective budget of uncertainty and thus our approach would lead to less conservative solutions than those of the conventional approach. 

}

\section{The IEEE Reliability Test System (RTS)} \label{rts24}
Figure~\ref{fig rts} illustrates the diagram of the IEEE RTS with 24 buses, 11 generators, and 4 wind farms. Loads are shown by arrows. Detailed data of the test system can be found in \cite{grigg1999ieee}.
\begin{figure}[h!] 
    \centering
    \includegraphics[width=0.5\textwidth , height=0.42\textheight
    ] {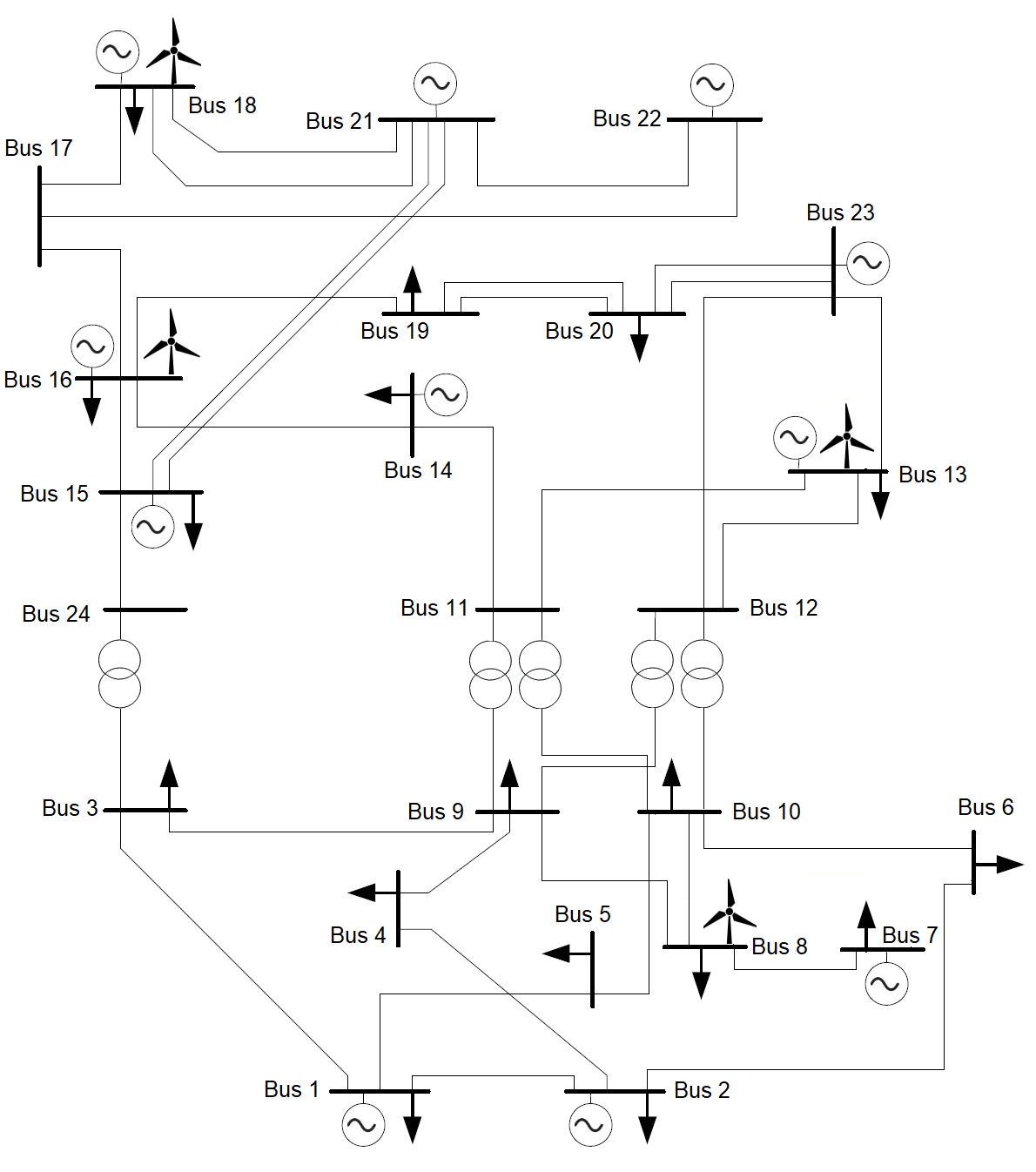}
    \caption{Diagram of the IEEE-RTS with 24 buses }
\label{fig rts}
\end{figure}
\end{APPENDICES}

\end{document}